\documentclass[letterpaper]{article}

\usepackage{arxiv}

\usepackage{style}

\usepackage{lineno,hyperref}

\usepackage{amsmath, amssymb, latexsym, multicol, amsthm, bm, placeins}
\usepackage[export]{adjustbox}
\usepackage{color}
\usepackage{multirow}

\usepackage{style}

\makeatletter\@addtoreset{equation}{section}\makeatother

\renewcommand{\shorttitle}{3D AVEM with stab-free a posteriori error bounds}

\numberwithin{equation}{section}

\makeatletter
\renewcommand{\@biblabel}[1]{#1\hfill \hspace{0.2cm}}
\makeatother

\begin{document}

\title{3D Adaptive VEM with stabilization-free\\ a posteriori error bounds}

\author{%
  Stefano Berrone\affil{1,}\corrauth
  and
  Davide Fassino\affil{1}
  and
  Fabio Vicini\affil{1}
}

\author{S. Berrone\thanks{Department of Mathematical Sciences “Giuseppe Luigi Lagrange”, Politecnico di Torino, Corso Duca degli Abruzzi 24, Torino, 10129, Italy. stefano.berrone@polito.it (S. Berrone), davide.fassino@polito.it (D. Fassino), fabio.vicini@polito.it (F. Vicini).}
	\And
	D. Fassino\footnotemark[1] \And F. Vicini\footnotemark[1]
}

\maketitle
\begin{abstract} 
The present paper extends the theory of Adaptive Virtual Element Methods (AVEMs) to the three-dimensional meshes showing the possibility to bound the stabilization term by the residual-type error estimator.
This new bound enables a stabilization-free a posteriori control for the energy error.
Following the recent studies for the bi-dimensional case, we investigate the case of tetrahedral elements with aligned edges and faces.
We believe that the AVEMs can be an efficient strategy to address the mesh conforming requirements of standard three-dimensional Adaptive Finite Element Methods (AFEMs), which typically extend the refinement procedure to non-marked mesh cells.
Indeed, numerical tests on the Fichera corner shape domain show that this method can reduce the number of three-dimensional cells generated in the refinement process by about $30\%$ with compared to standard AFEMs, for a given error threshold.
\end{abstract}


\maketitle
\section{Introduction}

Adaptive methods for approximating the solutions of partial differential equations are fundamental tools in scientific computing and engineering, whose mathematical background has been studied for decades \cite{CARSTENSEN20141195,NoSiVe:09, Bonito2024adaptive}. The standard adaptive methods are based on the iteration of the classical paradigm
\begin{equation*}
    {\sf SOLVE}\rightarrow{\sf ESTIMATE}\rightarrow{\sf MARK}\rightarrow{\sf REFINE},
\end{equation*}
to achieve at each loop a better approximation of the solution. The main idea behind the adaptivity resides in the definition of an a posteriori error estimator which indicates which elements bring the major error contribution, and then need to be refined. From the conformity of the mesh required by standard Adaptive Finite Element Methods (AFEMs), the refinement of the elements marked imposes the refinement of other elements just for geometric reasons. Since their introduction \cite{BasicPrinciples,TheHG} the Virtual Element Methods (VEMs), allowing the presence of aligned edges and faces, are a valuable strategy to allow the refinement of the marked element only. 
However, the extension of the adaptive FEM theory to the VEM is not trivial. 
The first contraction results for an adaptive VEM algorithm has been recently proposed in \cite{AVEMstabfree}  and the study of its optimality in \cite{AVEMConvergenceOptimality} in the lowest order case for triangular elements admitting the presence of aligned edges. 
The extension to the higher-order case has been introduced in \cite{higherorder}. 
One of the main difficulties in the study of an adaptive algorithm stands in the definition of a stabilization-free a posteriori error bound.
Indeed, the original VEM a posteriori error analysis showed in \cite{Cangiani} presented a bound of the type
\begin{equation}\label{eq:introErrorBound}
 \EnergyNorm{u-u_\mesh}\lesssim \eta_\mesh(u_\mesh) + S_\mesh(u_\mesh),   
\end{equation}
where $\EnergyNorm{u-u_\mesh}$ is the error in the energy norm between the solution $u$ of the continuous problem and the discrete solution $u_\mesh$ computed, $\eta_\mesh(u_\mesh)$ is residual type estimator and $S_\mesh(u_\mesh)$ is a stability term that needs to be added to the discrete bilinear form to get the coercivity. In particular, as remarked in \cite{StabMascotto}, this latter term presents two main complications: (i)  the stabilization term is not robust with respect to the polynomial degree of the method, and (ii) one should be really careful in the choice of $S_\mesh(u_\mesh)$ since the refinement procedure may produce elements with very small edges and faces. One of the first studies of an a posteriori error bound that does not involve the stability term, under specific assumptions, can be traced back to \cite{BerroneBorio2017}. Two of the possible strategies to overcome the presence of the stability term in \eqref{eq:introErrorBound} are the analysis of the recently introduced stabilization-free VEM schemes \cite{BerroneBorioMarconTeora,BerroneBorioMarcon} or the definition of stabilization-free a posteriori error bounds for the classical VEM. In this work we investigate the latter scenario. The bi-dimensional first-order VEM bound of the type 
\begin{equation*}
    S_\mesh(u_\mesh)\lesssim \eta_\mesh(u_\mesh),
\end{equation*}
is due to \cite{AVEMstabfree}. 
The purpose of this paper is to extend this result to the three-dimensional case and, consequently, to get 
\begin{equation*}
\EnergyNorm{u-u_\mesh}\lesssim \eta_\mesh(u_\mesh).
\end{equation*}
This stabilization-free a posteriori error bound brings to the definition of an adaptive algorithm, whose convergence result retraces the one of \cite{AVEMstabfree}.

The main novelty that arises in the higher-dimensional case is a new proof of the scaled Poincaré inequality, which results essential to prove the stabilization-free a posteriori error bound. In this paper we have decided to report the main propositions and theorems introduced in \cite{AVEMstabfree} for sake of readability, highlighting the differences which occur. Furthermore, we focus on the case of tetrahedral elements with aligned faces and aligned edges, since a refinement strategy for a general three-dimensional polyhedron is missing, to the best of our knowledge. 

The paper is structured as follows. In Section \ref{Sec:TheProblem} the continuous and the discrete problems are presented. We give here the definition of the global index of a node for the three dimensional case and the two preliminary properties on which the rest of the paper is based, i.e. the scaled Poincaré inequality and the Clément interpolation estimate. Section \ref{Sec:AposteroriErrorAnalysis} presents the stabilization-free a posteriori error bound and Section \ref{Sec:TheModuleGalerkin} discusses the convergence of the method. Finally, Section \ref{sec:NumericalExperiments} is devoted to the numerical experiments.
\section{The problem and the discretization}\label{Sec:TheProblem} We focus on the second-order \\Dirichlet boundary-value problem on a polyhedral domain $\Omega \subset \R^3$,
\begin{align}
    \begin{cases}
         -\nabla \cdot (\K \nabla u) + c u = f &\text{ in $\Omega$,}\\
         u = 0 &\text{ on $\partial \Omega$,}
    \end{cases}
    \label{Initial_Problem}
\end{align}
where the data, $\data = (\K, c, f)$, are $\K\in (L^\infty(\Omega))^{3\times 3}$ symmetric and positive-definite in $\Omega$, $c\in L^\infty(\Omega)$ non negative in $\Omega$, and $f\in L^2(\Omega)$. For simplicity, we consider homogeneous Dirichlet boundary conditions.
The extension to more general boundary conditions do not present theoretical difficulties. From problem \eqref{Initial_Problem}, we can derive its variational formulation, which reads as
\begin{align}
\begin{cases}
\text{find }u \in \mathbb{V}  \coloneqq H^1_{0}(\Omega) &\text{ such that}\\
\mathcal{B}(u,v) = (f,v), &\forall \;v \in \mathbb{V},
\end{cases}
\label{Variational_Problem}
\end{align}
where $(\cdot,\cdot)$ is the scalar product in $L^2(\Omega)$ and $\B(u,v) \coloneqq a(u,v) + m(u,v) $ is the bilinear form associated with problem, i.e.
  \begin{equation*}
    a(u,v) \coloneqq (\K\; \nabla u , \nabla v), \quad\quad \quad m(u,v) \coloneqq ( c\; u,v).
\end{equation*}  
From the assumption on the data, $\B$ is continuous and coercive, and problem \eqref{Variational_Problem} admits the existence and uniqueness of the solution $u$.
We indicate the energy norm as $\EnergyNorm{\cdot}\coloneqq \sqrt{\B(\cdot,\cdot)}$, which satisfies
\begin{equation}\label{bound_energynorm}
    c_\B \normH{v}^2\le \EnergyNorm{v}^2 \le c^\B\normH{v}^2, \quad \quad \quad\forall v \in \mathbb{V},
\end{equation}
where $\normH{\cdot}$ is the $H^1$-seminorm and $c_\B$ and $c^\B$ are two suitable strictly positive constants.

In the perspective of presenting a discrete formulation of the variational problem, we first introduce an initial conforming partition $\mesh_0$  of $\overline{\Omega}$ made of tetrahedra which satisfies the initial conditions given in \cite{StevensonNVB}.
Let $\mesh$ be any successive partition obtained from $\mesh_0$ by the refinement of some elements via the three-dimensional newest-vertex bisection (NVB) algorithm \cite{NoSiVe:09}. 

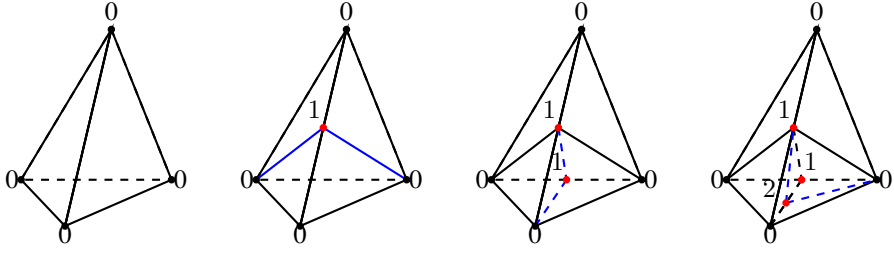
\begin{figure}[]
\centering
\begin{tikzpicture}[scale = 0.4]
\draw [thick] plot coordinates{  (2,0 ,0) (0,5 ,0) (0,0 ,4)  (2,0 ,0)} ;
\draw [thick] plot coordinates{ (-3, 0, 0) (0,5,0) (0,0,4) (-3, 0, 0)};
\draw [thick,dashed] plot coordinates{ (-3, 0, 0)  (0,5,0)  (2,0,0) (-3, 0, 0)} ;
\draw [thick] plot coordinates{ (0, 5, 0)  (0,0,4) } ;
\draw [very thick]plot [mark=*, only marks]coordinates{(-3,0,0)  (0,5,0) (0,0,4)  (2,0,0) };
\coordinate [label= 0] (z1) at (2.3,-0.6 ,0);
\coordinate [label= 0] (z1) at (-3.3,-0.6 ,0);
\coordinate [label= 0] (z1) at (0,5 ,0);
\coordinate [label= 0] (z1) at (0,-0.9,4);
\end{tikzpicture}
\hspace{0.3cm}
\begin{tikzpicture}[scale = 0.4]
\draw [thick] plot coordinates{ (2,0 ,0)  (0,5 ,0) (0,0 ,4)  (2,0 ,0)} ;
\draw [thick] plot coordinates{ (-3, 0, 0) (0,5,0) (0,0,4) (-3, 0, 0)};
\draw [thick,dashed] plot coordinates{ (-3, 0, 0)  (0,5,0)  (2,0,0) (-3, 0, 0)} ;
\draw [thick] plot coordinates{ (0, 5, 0)  (0,0,4) } ;
\draw [thick, blue] plot coordinates{(0, 2.5,2) (-3,0,0) } ;
\draw [thick, blue] plot coordinates{(0, 2.5,2) (2,0,0) } ;
\draw [very thick, red]plot [mark=*, only marks]coordinates{(0, 2.5,2)};
\draw [very thick]plot [mark=*, only marks]coordinates{(-3,0,0)  (0,5,0) (0,0,4)  (2,0,0) };
\coordinate [label= 0] (z1) at (2.3,-0.6 ,0);
\coordinate [label= 0] (z1) at (-3.3,-0.6 ,0);
\coordinate [label= 0] (z1) at (0,5 ,0);
\coordinate [label= 0] (z1) at (0,-0.9,4);
\coordinate [label= $1$] (z1) at (-0.3,2.5,2);
\end{tikzpicture}
\hspace{0.3cm}
\begin{tikzpicture}[scale = 0.4]
\draw [thick] plot coordinates{ (2,0 ,0)  (0,5 ,0) (0,0 ,4)  (2,0 ,0)} ;
\draw [thick] plot coordinates{ (-3, 0, 0) (0,5,0) (0,0,4) (-3, 0, 0)};
\draw [thick,dashed] plot coordinates{ (-3, 0, 0)  (0,5,0)  (2,0,0) (-3, 0, 0)} ;
\draw [thick] plot coordinates{ (0, 5, 0)  (0,0,4) } ;
\draw [thick] plot coordinates{(0, 2.5,2) (-3,0,0) } ;
\draw [thick] plot coordinates{(0, 2.5,2) (2,0,0) } ;
\draw [thick,blue,dashed] plot coordinates{(0, 2.5,2)  (-0.5,0,0) } ;
\draw [thick,blue,dashed] plot coordinates{(0, 0,4)  (-0.5,0,0) } ;
\draw [very thick, red]plot [mark=*, only marks]coordinates{(0, 2.5,2) (-0.5,0,0)};
\draw [very thick]plot [mark=*, only marks]coordinates{(-3,0,0)  (0,5,0) (0,0,4)  (2,0,0) };
\coordinate [label= 0] (z1) at (2.3,-0.6 ,0);
\coordinate [label= 0] (z1) at (-3.3,-0.6 ,0);
\coordinate [label= 0] (z1) at (0,5 ,0);
\coordinate [label= 0] (z1) at (0,-0.9,4);
\coordinate [label= $1$] (z1) at (-0.3,2.5,2);
\coordinate [label= $1$] (z1) at (-0.8,0,0);
\end{tikzpicture}
\hspace{0.3cm}
\begin{tikzpicture}[scale = 0.4]
\draw [thick] plot coordinates{ (2,0,0)  (0,5,0) (0,0,4)  (2,0,0)} ;
\draw [thick] plot coordinates{ (-3, 0, 0) (0,5,0) (0,0,4) (-3, 0, 0)};
\draw [thick,dashed] plot coordinates{ (-3, 0, 0)  (0,5,0)  (2,0,0) (-3, 0, 0)} ;
\draw [thick] plot coordinates{ (0, 5, 0)  (0,0,4) } ;
\draw [thick] plot coordinates{(0, 2.5,2) (-3,0,0) } ;
\draw [thick] plot coordinates{(0, 2.5,2) (2,0,0) } ;
\draw [thick,dashed] plot coordinates{(0, 2.5,2)  (-0.5,0,0) } ;
\draw [thick,dashed] plot coordinates{(0, 0,4)  (-0.5,0,0) } ;
\draw [thick,dashed] plot coordinates{(0, 0,4)  (-0.5,0,0) } ;
\draw [thick,dashed, blue] plot coordinates{(-0.25,0,2) (0, 2.5,2) } ;
\draw [thick,dashed, blue] plot coordinates{(-0.25,0,2) (2, 0,0) } ;
\draw [very thick, red]plot [mark=*, only marks]coordinates{(0, 2.5,2) (-0.5,0,0) (-0.25,0,2) };
\draw [very thick]plot [mark=*, only marks]coordinates{(-3,0,0)  (0,5,0) (0,0,4)  (2,0,0) };
\coordinate [label= 0] (z1) at (2.3,-0.6 ,0);
\coordinate [label= 0] (z1) at (-3.3,-0.6 ,0);
\coordinate [label= 0] (z1) at (0,5 ,0);
\coordinate [label= 0] (z1) at (0,-0.9,4);
\coordinate [label= $1$] (z1) at (-0.3,2.5,2);
\coordinate [label= $1$] (z1) at (-0.2,0,0);
\coordinate [label= $2$] (z1) at (-0.8,-0.1,2) ;
\end{tikzpicture}
\caption{Three successive refinements of a tetrahedron. Blue lines represent the edges added at each refinement. Black points and red points represent the proper and the hanging nodes, respectively. The numbers are the global indices.}
\label{Figure:GlobalIndex}
\end{figure}

Given a tetrahedron $E\in \mesh_0$, this refinement consists in selecting an edge $\tilde{e}$ of $E$ and connecting its midpoint with the two remaining vertices ``opposite'' to $\tilde{e}$. 
The two new edges so built are part of the new triangular face formed in the refinement procedure that splits the tetrahedron into two elements (see Figure \ref{Figure:GlobalIndex} for a graphical representation). 

We remark that we have decided to use the term {\it tetrahedron} to denote the three-dimensional simplex, and we denote by {\it tetrahedral-shape elements} all the elements obtained during the refinement process that might contain aligned edges and faces.

Given a tetrahedral-shape element $E \in \mesh$, we define a {\textit {straight side}} of $E$ each one of the six edges of the smallest simplex that contains $E$.
Moreover, let $\nodes_E$ denote the collection of all the nodes of $E$.
We refer to each of the four endpoints of the {\textit {straight side}} of $E$ as {\textit {proper node of E}}, or simply {\textit {vertex}}.
These four vertices are collected in the set $\propers_E$, namely the set of {\textit {proper nodes of $E$}}.
Finally, we indicate by $\hangings_E$ the set of all the nodes sitting on the boundary of $E$ that are not {\textit {propers}} for $E$, i.e., $\hangings_E := \nodes_E \setminus \propers_E$.
Note that, in the classic case where $E$ is a tetrahedron, the vertices and the nodes of $E$ coincide. 
Specifically, $\nodes_E = \propers_E$ and $\hangings_E = \emptyset$.

At a global level, we denote by $\nodes$ the set of all the nodes of the partition, we define the \textit{proper} nodes set as
$\mathcal{P}:= \left\{\bm{x}\in \mathcal{N}: \bm{x}\in \mathcal{P}_E, \ \forall E \text{ containing }\bm{x}\right\}$, and we call \textit{hanging} the nodes that are not \textit{proper} nodes, i.e. $\hangings = \nodes \setminus \propers$. 
We recall that the presence of hanging nodes at the boundary of the elements $E$, can be handled by the VEM theory, which treats tetrahedra with hanging nodes as polyhedra with aligned faces and edges. 
 We denote by $\mathcal{F}_E$ the set of faces of each element $E \in \mesh$, by $\partial F$ the boundary of the face $F$, and by $\mathcal{E}_F$ the set of edges of any face $F$. In the paper, we set $h_E = |E|^{1/3}$, where $|E|$ denotes the volume of $E$, and $h=\max_{E\in \mesh} h_E$. 

A crucial concept, firstly introduced in \cite{AVEMstabfree} for the bi-dimensional case and extended in \cite{higherorder}, is the \textit{global index of a node}. To define it, we recall that in the NVB refinement strategy, a hanging node $\x \in \hangings$  is always obtained by bisecting an edge of the tetrahedron. We denote  by $\{\bm{x'}, \bm{x''}\}$ the 
 set of the endpoints of the refined edge for which $\x$ is the midpoint.

\begin{defn}[global index of a node]\label{def:global_Index} Given a node $\x \in \nodes$, we define its global index $\lambda$ as:
\begin{itemize}
    \item $\lambda(\x)=0,$ if $\x \in \propers$;
    \item $\lambda(\x)= \max\{\lambda(\bm{x'}), \lambda(\bm{x''})\}+1$, if $\x \in \hangings$.
\end{itemize}\end{defn}
An example of the evolution of the global index of a node, after three refinements is shown in Figure \ref{Figure:GlobalIndex}.

In the analysis presented, we suppose that the largest global index of the partition does not blow. More specifically, we have the following assumption.
\begin{assump}\label{ass:LambdaMax}
    We assume the existence of a constant $\Lambda\ge1$ such that 
\begin{equation*}\label{eq:Lambda-admission}
    \Lambda_\mesh\coloneqq \max_{\x \in \nodes}\lambda(\x)\le \Lambda.
\end{equation*}
We call this type of partition $\Lambda$-admissible.
\end{assump}

\begin{remark}\label{Remark:LambdaMax} The case $\Lambda=0$ is not included in Assumption~\ref{ass:LambdaMax} and it is not considered in the following analysis. Indeed, if $\Lambda=0$ the partition does not admit hanging nodes, then the VEM scheme reduces to the FEM one.
\end{remark}
\begin{remark}\label{r:condeguences_of_lambda}
    If a partition $\mesh$ is $\Lambda$-admissible, the number of hanging nodes on each side of an element can be upper bounded by $2^\Lambda-1$ hanging nodes (see also \cite[Remark 2.3]{AVEMstabfree}).
\end{remark}
\subsection{VEM spaces and projectors}
The construction of the three-dimensional VEM spaces is an extension of 
the bi-dimensional ones \cite{EquivalentProjectors,  BasicPrinciples, TheHG}.  
Let $E$ be a tetrahedral-shape element in the partition $\mesh$, and $F$ a face in $\mathcal{F}_E$. On $F$, we define the enhanced bi-dimensional space of order one, starting from the boundary of $F$
\begin{equation*}
    \V_{\partial F}\coloneqq \{ v \in C^0(\partial F): v|_e \in \P_1(e),\;  \forall e\in \mathcal{E}_F\},
\end{equation*}
and extending it in the interior of $F$,
\begin{equation*}
\begin{split}
    \V_F \coloneqq\{v \in H^1(F): v|_{\partial F} \in \V_{\partial F}, \;\Delta v|_F \in \P_1(F),\;
    \\  \int_F \left(v- \Pi^\nabla_F v\right)q_1=0 ,\; \forall q_1 \in \P_1(F)\},
    \end{split}
\end{equation*}
where the projector $\Pi^\nabla_F : H^1(F)\rightarrow \P_1(F)$ is such that $\forall v \in H^1(F)$
\begin{equation*}
    \int_F \nabla \left( \Pi^\nabla_F v  - v \right)\nabla q_1 = 0\quad \forall q_1 \in \P_1(F),\quad \quad \int_{\partial F}\left(\Pi^\nabla_F v -v\right) =0.
\end{equation*}
In the element $E$, we define the space $\V_E$ such that it reduces to $\V_F$ on each face of $\mathcal{F}_E$. In particular, we define on the boundary of $E$,
\begin{equation*}
\V_{\partial E}\coloneqq \{v \in C^0(\partial E): \; v|_F \in \V_F,\; \forall F \in\mathcal{F}_E\},
\end{equation*}
and the enhanced local space 
\begin{equation}\label{def:V_E}
\begin{split}
    \V_E \coloneqq\{v \in H^1(E): \; v|_{\partial E}\in \V_{\partial E}, \;\Delta v|_E \in \P_1(E), \; \\\int_E \left( v - \Pi^\nabla_ E\right) p_1 =0, \; \forall p_1 \in \P_1(E)\},\end{split}
\end{equation}
in which $\Pi^\nabla_E: H^1(E)\rightarrow \P_1 (E)$ is defined as $\forall v \in H^1(E)$
\begin{equation}\label{cond_PiNabla1}
    \int_E \nabla \left( \Pi^\nabla_E v  - v \right)\nabla w = 0\quad \forall w \in \P_1(E), \quad \quad\quad  \int_{\partial E}\left(\Pi^\nabla_E v -v\right) =0.
\end{equation}
We denote the global space as
\begin{equation*}
    \Vmesh \coloneqq\{v \in \V: v|_E \in \V_E,\; \forall E \in \mesh\},
\end{equation*}
the space of piecewise polynomials of order one on $\mesh$ as
\begin{equation*}
    \Wmesh\coloneqq\{w \in L^2(\Omega): w|_E\in \P_1 (E),\; \forall E\in \mesh\},
\end{equation*}
and the subspace $\Vmeshz$, i.e.
\begin{equation}\label{def:V_meshz}
  \Vmeshz \coloneqq \Vmesh \cap \Wmesh,  
\end{equation}
which plays a crucial role in the proofs. 
\begin{remark}
    Notice that in general  $\Wmesh\not \subset \Vmesh $,  indeed $ \Vmesh \subset C^0(\Omega)$, while if $v\in\Wmesh$ in general $v\not \in C^0(\Omega)$, but only piece-wise polynomial on each element $E$ of the partition $\mesh$. 
\end{remark}
Let us define the projection operator $\Pi^\nabla_\mesh: \Vmesh \rightarrow \Wmesh$ such that when we restrict it to an element $E\in\mesh$ we get $\Pi^\nabla_E$. Finally, let $\Ical_E:\V_E\rightarrow \P_1(E)$ be the Lagrange interpolation operator at the nodes of $E$ and the respective global Lagrange interpolation operator $\Icalmesh:\Vmesh \rightarrow \Wmesh$ such that $\Icalmesh|_E = \Ical_E$. 

We briefly recall some properties of the projection operator $\Pi^\nabla_E$, which will be used in the paper.
\begin{lemma}
For the projection operator $\Pi^\nabla_E$ it holds the following
\begin{equation}\label{eq:continuity_PiNabla}
        \normLE{\Pi^\nabla_E v}\lesssim \normLE{v} + h_E\normHE{v},\quad \quad \forall v \in \V
\end{equation}
where the hidden constant does not depend on the diameter $h_E$. Furthermore, 
\begin{equation}\label{diseq:PiNabla}
    \normHE{v- \Pi^\nabla_E v}\le \normHE{v-\Ical_E v}, \quad\quad \forall v \in \V.
\end{equation}
\end{lemma}
\begin{proof}
   
    The proof of \eqref{eq:continuity_PiNabla} can be found in \cite[Lemma 3.4]{BrennerSung}. On the other hand, \eqref{diseq:PiNabla} is a direct consequence of the  $H^1$ orthogonality of operator $\Pi^\nabla_E$ \eqref{cond_PiNabla1}, i.e.
    \begin{align*}
        \normHE{v- \Pi^\nabla_E v}^2 &=  \left( \nabla\left(v- \Pi^\nabla_E v\right)_E, \nabla\left(v- \Pi^\nabla_E v\right)\right)_E
        =\left( \nabla\left(v- \Pi^\nabla_E v\right), \nabla v\right)_E
        \\&= \left( \nabla\left(v- \Pi^\nabla_E v\right), \nabla (v - \Ical_E v)\right)_E\le \normHE{v- \Pi^\nabla_E v}\normHE{v- \Ical_E v}.
    \end{align*}
\end{proof}

 \subsection{The discrete problem}
We assume the data to be piecewise constant on each element $E$ of the partition. In the following, we denote the restriction of the data $(\K, f,c)$ over the element $E$ as $(\K_E, f_E, c_E)$.
We now define the discrete bilinear forms $a_\mesh$ and $m_\mesh:$ $ \Vmesh\times \Vmesh \rightarrow \R$, such that
\begin{align*}
    &a_E(v,w)\coloneqq\int_E \K_E \left(\nabla \Pi^\nabla_E v\right) \cdot \left(\nabla \Pi^\nabla_E w\right), &a_\mesh (v,w) \coloneqq  \sum_{E\in \mesh}a_E(v,w); \\
    &m_E(v,w)\coloneqq\int_E c_E \ \Pi^\nabla_E v \;\Pi^\nabla_E w, &m_\mesh (v,w) \coloneqq  \sum_{E\in \mesh}m_E(v,w).
\end{align*}
We introduce the local stabilization form $s_E: \V_E\times \V_E \rightarrow \R$ as
\begin{equation}\label{eq:StabilizationE}
    s_E (v,w) = h_E\sum^{\nodes_E}_{i=1} v(\x_i) w(\x_i),
\end{equation}
where $\{\x_i\}_{i=1}^{\mathcal{N}_E}$ is the set of the degrees of freedom of $E$.
We remark that this classical choice of the stabilization term differs from the bi-dimensional classical one for the presence of the $h_E$ factor. Furthermore, we recall that the choice of the form of the stabilization term is not restrictive and the results presented can be extended to other types of the stabilization term. In \cite{BrennerSung} it is shown that, for the choice \eqref{eq:StabilizationE}, it holds
\begin{equation}\label{e:stab_equivalence}
 c_s \normHE{v}^2\le s_E(v,v) \le C_s\normHE{v}^2,   \quad \forall v\in \V_E/ \mathbb{R},
\end{equation}
with $C_s$ and $c_s$ positive constants. We define the stabilization form on the partition $\mesh$ as
\begin{equation*}
    S_\mesh(v,w)\coloneqq \sum_{E\in \mesh}s_E(v-\Ical_E v, w - \Ical_E w)\qquad\forall v,w\in\Vmesh,
\end{equation*}
which yields to
\begin{equation}\label{eq:S_mesh_Hnorm}
    S_\mesh (v,v)\simeq \normHT{v - \Icalmesh v}^2\qquad\forall v\in \Vmesh,
\end{equation}
where $\normHT{\cdot}$ is the broken $H^1$-seminorm over the mesh $\mesh$.
We have now all the tools to define the complete bilinear form $\B_\mesh:\Vmesh\times \Vmesh \rightarrow \R$, associated to problem \eqref{Initial_Problem},
\begin{equation}\label{def:definition_BT}
    \B_\mesh(v,w)\coloneqq a_\mesh (v,w) + m_\mesh (v,w) + \gamma S_\mesh(v,w), 
\end{equation}
where $\gamma \ge \gamma_0 > 0$. 
Finally, to approximate the forcing term we define $\Fcalmesh: \Vmesh \rightarrow \R$ such that
\begin{equation}\label{def:definition_FT}
    \Fcalmesh(v)\coloneqq \sum_{E\in \mesh}\int_E f_E \Pi^\nabla_E v, \qquad \forall\,v\, \in \Vmesh.
\end{equation}
We define the discrete formulation of the problem as 
\begin{align}
\begin{cases}
\text{find }u_\mesh \in \Vmesh &\text{such that}\\
\B_\mesh(u_\mesh, v)=\Fcalmesh(v), &\forall \;v \in \Vmesh,
\end{cases}
\label{Discrete_Problem}
\end{align}
which admits a unique solution.

The importance of $\Vmeshz$ can be already viewed in the following Lemma, whose proof can also be found in \cite[Lemma 2.6]{AVEMstabfree}.
\begin{lemma}[Galerkin quasi-orthogonality] \label{lemma:Galerkin_Quasi_Ort}
    For any $v\in \Vmeshz$, it holds
    \begin{equation*}
        \B(u - u_\mesh, v)=0,
    \end{equation*}
    where $u$ is the solution of \eqref{Variational_Problem} and $u_\mesh$ the solution of the discrete problem \eqref{Discrete_Problem}.
\end{lemma}
\begin{proof} From the continuous and the discrete formulation of the variational problem we get, $\forall v \in \Vmesh^0$,
\begin{equation*}
    \B(u- u_\mesh,v)= \B(u, v)- \B(u_\mesh, v)= (f,v)- \Fcalmesh(v)+\B_\mesh (u_\mesh, v)-\B (u_\mesh,v).
\end{equation*}
Since $v\in\Vmesh^0$, $\Pi^\nabla_\mesh v = v$, we have from \eqref{def:definition_FT} $\Fcalmesh(v)= (f,v)$. Furthermore, we notice that from \eqref{cond_PiNabla1} we get
\begin{align*}
    a_\mesh(u_\mesh, v) = \sum_{E\in \mesh}\int_E \K_E\left (\nabla \Pi^\nabla_E u_\mesh \right ) \cdot\left (\nabla \Pi^\nabla_E v \right ) &=\sum_{E\in \mesh}\int_E \K_E\left (\nabla \Pi^\nabla_E u_\mesh \right ) \cdot\left (\nabla v \right ) \\&= \sum_{E\in \mesh}\int_E \K_E\left (\nabla u_\mesh \right ) \cdot\left (\nabla v \right )= a(u_\mesh,v).
\end{align*}
Finally, from the definition of the enhanced Virtual Element space \eqref{def:V_E}, $\forall v \in \Vmeshz$
\begin{align*}
    m_\mesh(u_\mesh, v) = \sum_{E\in \mesh}\int_E c_E \Pi^\nabla_E u_\mesh  \;\Pi^\nabla_E v &=\sum_{E\in \mesh}\int_E c_E \Pi^\nabla_E u_\mesh\;   v \\&= \sum_{E\in \mesh}\int_E c_E u_\mesh \;v = m(u_\mesh,v),
\end{align*}
which concludes the proof.
\end{proof}
\subsection{The preliminary properties}\label{sec:TheKeyProperties}
We present two key properties that will be essential in the rest of the paper. The first of these properties is the scaled Poincaré inequality, which retraces the result presented in the bi-dimensional case \cite[Proposition 3.1]{AVEMstabfree}, but with crucial differences in the proof. We present previously the following Lemmas.

\begin{lemma} \label{lemma:NotTheSameLambda}
Let $E$ be an arbitrary element in $\mesh$, $E$ cannot have all its vertices with the same global index $\Tilde{\lambda}>0$.
\end{lemma}

\begin{figure}[]
        \centering        
         \begin{tikzpicture}[scale = 0.7]
         \fill[bubbles] plot coordinates{ (-3, 0, 0)  (0, 2.5,2)  (2,0,0) (-3, 0, 0)} ;
         \fill[bubbles] plot coordinates{ (-3, 0, 0)  (0,0,4)  (2,0,0) (-3, 0, 0)} ;       
         \draw [thick] plot coordinates{ (2,0,0)  (0,5,0) (0,0,4)  (2,0,0)} ;
         \draw [thick] plot coordinates{ (-3, 0, 0) (0,5,0) (0,0,4) (-3, 0, 0)};
         \draw [thick,dashed] plot coordinates{ (-3, 0, 0)  (0,5,0)  (2,0,0) (-3, 0, 0)} ;
         \draw [thick] plot coordinates{ (0, 5, 0)  (0,0,4) } ;
         \draw [thick] plot coordinates{(0, 2.5,2) (-3,0,0) } ;
         \draw [thick] plot coordinates{(0, 2.5,2) (2,0,0) } ;
         \draw [very thick]plot [mark=*, only marks]coordinates{(0, 2.5,2) (2,0,0)  (0, 5, 0)  (0,0,4) (-3, 0, 0)};
         \draw [very thick,red]plot [mark=*, only marks]coordinates{(0, 2.5,2) };
         \coordinate [label= $\lambda_j$] (z1) at (2.4,-0.6 ,0);
         \coordinate [label= $\lambda_m$] (z1) at (-3.2,-0.6 ,0);
         \coordinate [label= $\lambda_\ell$] (z1) at (0.5,4.9 ,0);
          \coordinate [label= $\x'_0$] (z1) at (-0.2,5 ,0);
         \coordinate [label= $\lambda_k$] (z1) at (0.4,-0.7,4);
         \coordinate [label= $\x''_0$] (z1) at (-0.4,-0.8,4);
         \coordinate [label= $\lambda_i$,red] (z1) at (0.4,2.5,2);
         \coordinate [label= $\bm{x}_0$,red] (z1) at (-0.3,2.5,2);
         \coordinate [label= $T$] (z1) at (0.9,2.9,2);
         \coordinate [label= $E$] (z1) at (0.4,1.1,2);
    \end{tikzpicture}
    \caption{Element $E\in \mesh$ obtained after the refinement of element $T$ via the NVB, where the edge bisected is the one with endpoints $\x'_0$ and $\x''_0$. $\lambda_j,\lambda_k,\lambda_\ell,\lambda_m$ are the global indices of the vertices of the element $T$, whereas $\lambda_i,\lambda_j,\lambda_k,\lambda_m$ are the global indices of the vertices of $E$.}
     \label{Figure:Hangings}
     \end{figure}
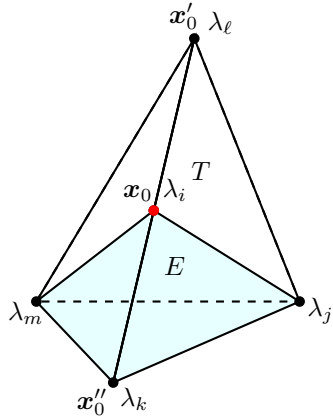

\begin{proof}
     We prove it by contradiction. By hypothesis, $E$ has only hanging nodes as vertices ($\Tilde{\lambda}>0$), then it is generated after a partition. We reference to Figure~\ref{Figure:Hangings}, one of the vertices of $E$, say $\x_{0}$, is the midpoint of the edge with endpoints $\x'_{0}$ and $\x''_{0}$. Since one of these endpoints, say $\x''_0$, is another vertex of $E$, the global index of $\x''_0$ should be $\Tilde{\lambda}$. However, by Definition~\ref{def:global_Index}, the global index of $\x_{0}$ results
        \begin{equation*}
        \Tilde{\lambda}=\lambda(\x_{0}) = \max\{\lambda(\x'_{0}), \lambda(\x''_{0})\}+ 1\ge  \Tilde{\lambda} +1,
        \end{equation*}
        which yields to a contradiction.
\end{proof}

In particular, we notice that if we set $\Lambda =1$ in Assumption~\ref{ass:LambdaMax} all the elements of the partition contain at least a proper node.  

\begin{lemma} \label{lemma:reduceLambda}
Let $E$ be an element with only hanging nodes as vertices, that is generated by the refinement of an element $T$. Passing from $E$ to its parent $T$, the sum of the global indices of the vertices of $E$ strictly increases by at least one unit compared to the sum of the global indices of the vertices of $T$.
\end{lemma}
\begin{proof}
        In this proof, we refer again to Figure~\ref{Figure:Hangings}.
        We denote by $\lambda_i, \lambda_j, \lambda_k, \lambda_m$ the global indices of the vertices $i,j,k,m$ of $E$ respectively that are assumed to be strictly positive. Let $T$ be the parent of $E$ and $\lambda_\ell \ge 0$ the global index of the vertex $\ell$ of $T$ not belonging to $E$. We claim that 
        \begin{equation*}
        \lambda_\ell < \lambda_i.
        \end{equation*}
        To prove this, we observe that $i$ is the midpoint of the edge whose global indices are $\lambda_\ell$ and $\lambda_k$ (see Figure \ref{Figure:Hangings}). By Definition \ref{def:global_Index}, $\lambda_i= \max\{\lambda_\ell, \lambda_k\} +1$, whence if $\lambda_k\le \lambda_\ell$, then $\lambda_i =\lambda_\ell +1 > \lambda_\ell$, whereas if $\lambda_\ell <\lambda_k$, then $\lambda_i = \lambda_k+1 >\lambda_\ell +1 >\lambda_\ell$. 
\end{proof}

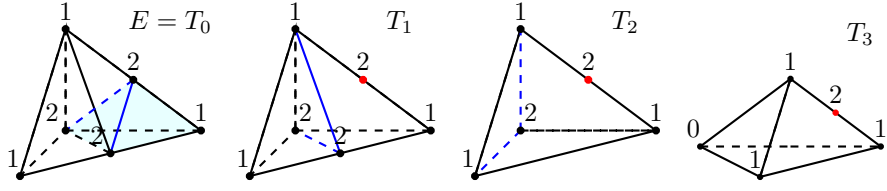
\begin{figure}[]
\centering
\begin{tikzpicture}[scale = 0.45]
\fill[bubbles]plot coordinates{(0,0 ,0)  (2,1.5,0)(2,0,1.75) (0,0 ,0) };
\fill[bubbles]plot coordinates{(4,0 ,0) (2,1.5,0)(2,0,1.75) (4,0 ,0)};
\draw [thick,dashed] plot coordinates{  (0,0 ,0) (4,0 ,0) (0,3 ,0)  (0,0 ,0)} ;
\draw [thick,dashed] plot coordinates{ (0, 0, 0) (0,3,0) (0,0,3.5) (0, 0, 0)};
\draw [thick] plot coordinates{ (4, 0, 0)  (0,3,0)  (0,0,3.5) (4, 0, 0)} ;
\draw [thick,dashed] plot coordinates{(2, 0,1.75) (0,0 ,0) } ;
\draw [thick] plot coordinates{(2, 0,1.75) (0,3 ,0) } ;
\draw [thick,dashed, blue] plot coordinates{((2,1.5,0) (0,0 ,0) } ;
\draw [thick,blue] plot coordinates{(2,1.5,0)(2,0,1.75) } ;
\draw [very thick]plot [mark=*, only marks]coordinates{(0,0 ,0) (4,0 ,0) (0,3 ,0) (0,0,3.5) (2,0,1.75) (2, 1.5,0) };
\coordinate [label= $2$] (z1) at (0-0.4,0 ,0);
\coordinate [label= $1$] (z1) at (4,0 ,0);
\coordinate [label= $1$] (z1) at (0,3 ,0);
\coordinate [label= $1$] (z1) at (-0.2,0,3.5);
\coordinate [label= $2$] (z1) at (2-0.4,0,1.75);
\coordinate [label= $2$] (z1) at (2,1.5,0);
\coordinate [label= {$E =T_0$}] (z1) at (4.25,3.75,3);
\end{tikzpicture}
\begin{tikzpicture}[scale = 0.45]
\draw [thick,dashed] plot coordinates{  (0,0 ,0) (4,0 ,0) (0,3 ,0)  (0,0 ,0)} ;
\draw [thick,dashed] plot coordinates{ (0, 0, 0) (0,3,0) (0,0,3.5) (0, 0, 0)};
\draw [thick] plot coordinates{ (4, 0, 0)  (0,3,0)  (0,0,3.5) (4, 0, 0)} ;
\draw [thick,dashed, blue] plot coordinates{(2, 0,1.75) (0,0 ,0) } ;
\draw [thick, blue] plot coordinates{(2, 0,1.75) (0,3 ,0) } ;
\draw [very thick,red]plot [mark=*, only marks]coordinates{(2, 1.5,0) };
\draw [very thick]plot [mark=*, only marks]coordinates{(0,0 ,0) (4,0 ,0) (0,3 ,0) (0,0,3.5) (2,0,1.75)};
\coordinate [label= $2$] (z1) at (0.3,0 ,0);
\coordinate [label= $1$] (z1) at (4,0 ,0);
\coordinate [label= $1$] (z1) at (0,3 ,0);
\coordinate [label= $1$] (z1) at (-0.2,0,3.5);
\coordinate [label= $2$] (z1) at (2,0,1.75);
\coordinate [label= $2$] (z1) at (2,1.5,0);
\coordinate [label= {$T_1$}] (z1) at (4.25,3.75,3);
\end{tikzpicture}
\begin{tikzpicture}[scale = 0.45]
\draw [thick,dashed] plot coordinates{  (0,0 ,0) (4,0 ,0) (0,0 ,0)} ;
\draw [thick,dashed] plot coordinates{ (0,3,0) (0,0,3.5) };
\draw [thick] plot coordinates{ (4, 0, 0)  (0,3,0)  (0,0,3.5) (4, 0, 0)} ;
\draw [thick,dashed, blue] plot coordinates{  (0,0 ,0) (0,3 ,0) } ;
\draw [thick,dashed, blue] plot coordinates{  (0,0 ,0) (0,0,3.5) } ;
\draw [very thick,red]plot [mark=*, only marks]coordinates{(2, 1.5,0) };
\draw [very thick]plot [mark=*, only marks]coordinates{(0,0 ,0) (4,0 ,0) (0,3 ,0) (0,0,3.5)};
\coordinate [label= $2$] (z1) at (0.3,0 ,0);
\coordinate [label= $1$] (z1) at (4,0 ,0);
\coordinate [label= $1$] (z1) at (0,3 ,0);
\coordinate [label= $1$] (z1) at (-0.2,0,3.5);
\coordinate [label= $2$] (z1) at (2,1.5,0);
\coordinate [label= {$T_2$}] (z1) at (4.25,3.75,3);
\end{tikzpicture}
\begin{tikzpicture}[scale = 0.3]
\draw [thick] plot coordinates{  (-4,0 ,0)(0,0,3.5) } ;
\draw [thick] plot coordinates{  (-4,0 ,0)(0,3,0) } ;
\draw [thick,dashed] plot coordinates{  (-4,0 ,0) (0,0 ,0)} ;
\draw [thick,dashed] plot coordinates{  (0,0 ,0) (4,0 ,0)} ;
\draw [thick,dashed] plot coordinates{ (0,3,0) (0,0,3.5) };
\draw [thick] plot coordinates{ (4, 0, 0)  (0,3,0)  (0,0,3.5) (4, 0, 0)} ;
\draw [very thick,red]plot [mark=*, only marks]coordinates{(2, 1.5,0) };
\draw [very thick]plot [mark=*, only marks]coordinates{(-4,0 ,0) (4,0 ,0) (0,3 ,0) (0,0,3.5)};
\coordinate [label= $0$] (z1) at (-4-0.3,0 ,0);
\coordinate [label= $1$] (z1) at (4,0 ,0);
\coordinate [label= $1$] (z1) at (0,3 ,0);
\coordinate [label= $1$] (z1) at (-0.2,0,3.5);
\coordinate [label= $2$] (z1) at (2,1.5,0);
\coordinate [label= {$T_3$}] (z1) at (4.25,5.25,3);
\end{tikzpicture}

\caption{ The four elements of the chain \eqref{eq:chain} from the light blue tetrahedral-shape element $E$ to the element with a proper node $T_3$. This chain can also be seen as three successive refinements of element $T_3$. Blue lines represent the edges added at each refinement. Black points represent the vertices of the considered elements, red points denote the hanging nodes. The numbers next to the vertices are their global indices. The fixed maximum global index is $\Lambda =2$.}
\label{Figure:Step4}
\end{figure}

\begin{lemma}\label{lemma:chainPoincare} Let $E$ be an element with only hanging nodes as vertices and let us denote the chain of ancestors of $E$ as:
        \begin{equation}\label{eq:chain}
           \Acal_N(E) =\{T_0 =E, T_1, \dots, T_N\}, 
        \end{equation}
        where $T_N$ is the first element containing a proper node.
    We claim that the length of this chain can be bounded by $3 (\Lambda -1)$, and it holds 
\begin{equation}\label{eq:rel_hT_hE}
         h_{T_N} \lesssim 2^{\Lambda -1}h_{E}.
     \end{equation}
    
\end{lemma}
\begin{proof}
    Let $E$ be an element fixed with only hanging nodes as vertices and $\Acal_N(E)$ the chain of its ancestors \eqref{eq:chain}.
    The number of elements $N$ of $\Acal_N(E)$ can be bounded as
    \begin{equation}\label{eq:boundLegnthChain}
    N\le 3(\Lambda -1). 
    \end{equation}
    With the intent to prove this, let us firstly consider the longest not-admissible chain (according to Lemma \ref{lemma:NotTheSameLambda})  starting from an element with vertices having the same global index $\Lambda$ to the one with the same global index $1$.  From Lemma~\ref{lemma:reduceLambda} this chain  has maximum length $4 (\Lambda -1)$. To ensure the admissiblity of the chain, according to Lemma~\ref{lemma:NotTheSameLambda}, we need to exclude all the $\Lambda$ elements of the chain with the same global index. In this way, we have obtained a chain that ends with an element with no proper nodes. Finally, to close the definition of  $\Acal_N(E)$ we add the element with at least one proper node, $T_N$, proving the inequality, i.e.
    \begin{equation*}
    N\le 4(\Lambda -1) -\Lambda +1 = 3(\Lambda-1).
    \end{equation*}
    See for instance Figure~\ref{Figure:Step4}, representing $\Acal_N(E)$ of the light-blue tetrahedral-shape element $E$ in a mesh in which $\Lambda =2$ and, from \eqref{eq:boundLegnthChain}, $N\le 3$. In the picture, $E$ needs three steps to reach the element $T_3$, i.e. the one having at least one proper node as a vertex. This chain is the longest admissible one, since a further refinement of $E$ would bring a global index equal to 3. Furthermore, we remark that the presence of hanging nodes in the elements of the chain does not influence on the length of the chain: as shown in Figure \ref{Figure:Step4}, the length of the chain from $T_2$ to $T_3$ does not consider the red point.

    Each time an element is refined, its volume halves, i.e. if $T$ is the parent of $E$, $h^3_{T} \simeq 2 h_E^3 $. Then, employing the bound \eqref{eq:boundLegnthChain}, for the smallest element of the chain \eqref{eq:chain} and the largest one it holds 
    \begin{equation}\label{eq:relazioneHcubo}
        h^3_{T_N} \lesssim 2^{3(\Lambda -1)}h^3_{E}.
    \end{equation}
    Taking the cube root, we get \eqref{eq:rel_hT_hE}.
\end{proof}

\begin{prop}[scaled Poincaré inequality in $\mathbb{V}_\mathcal{T}$]\label{ScaledPoincare}
    There exists a constant $C_P>0$ depending on $\Lambda$ and independent of $h_E$, such that    
    \begin{linenomath}
    \begin{align}\label{eq:ScaledPoicare}
     &\sum_{E\in \mathcal{T}}h_E^{-2}  \normLE{v}^2 \le C_P \normH{v}^2 &\forall v \in \mathbb{V}_\mathcal{T}\text{ such that } v(\bm{x})=0,\; \forall \x\in \mathcal{P}.
    \end{align}
    \end{linenomath}
\end{prop} 
\begin{proof} 

We assume $v\in \Vmesh$ such that $v(\x)=0$, $\forall \x \in \propers$, and we divide the proof into the following steps.
        Let $E\in \mesh$ be an arbitrary element. If at least one of its vertices is a proper node, we can use the classical Poincaré inequality \cite{Gilbard2001, NoSiVe:09}
        \begin{equation}\label{e:ScaledPoincare1}
             h_E^{-2}  \normLE{v}^2 \lesssim \normHE{v}^2.
        \end{equation}  
        On the other hand, if $E$ has all hanging nodes as vertices, as shown for instance in Figure \ref{Figure:Example_all_hangings}, we build for $E$ the chain $\Acal_{N}(E)$ \eqref{eq:chain} in Lemma \ref{lemma:chainPoincare}. Let $T_N$ be the last element of the chain, from \eqref{eq:rel_hT_hE} we get 
    \begin{equation*}
    h^{-2}_E\norm{v}^2_{0,E}\le h_E^{-2}\norm{v}^2_{0, T_N}= \left(\frac{h_{T_N}}{h_E}\right)^2 h^{-2}_{T_N} 
 \norm{v}^2_{0,T_N}\;  \lesssim 2^{2(\Lambda -1)}\abs{v}^2_{1,T_N},
    \end{equation*}
    in which the last inequality exploits the presence of a proper node in $T_N$.
    
    We introduce now $\TT_{T_N}$ the union set of all the elements $E$ generated by the refinements of $T_N$ (all these elements have a chain $\Acal_{N}(E)$ ending with $T_N$). 
    Recalling the inequality~\eqref{eq:relazioneHcubo}, the number of the elements of $\TT_{T_N}$ is bounded by $2^{3(\Lambda-1)}$.
    Thus, we have 
    \begin{equation*}
    \sum_{E \in \TT_{T_N}} h^{-2}_{E}\norm{v}^2_{0,E}\lesssim 2^{5(\Lambda-1)}\abs{v}^2_{1,T_N}.
    \end{equation*}
    Finally, summing up all the elements in $\mesh$, and denoting by $\mesh^\hangings$ the set of elements with only hanging nodes, we have
    \begin{equation*}
    \sum_{E\in \mesh} h^{-2}_E \norm{v}_{0,E} \lesssim \sum_{E\in \mesh\setminus\mesh^\hangings} |v|^2_{1,E} + \sum_{E\in \mesh^\hangings}h^{-2}_E \norm{v}^2_{0,E}\lesssim 2^{5(\Lambda -1)}\normH{v}^2,
    \end{equation*}
    which concludes the proof.
\end{proof}
      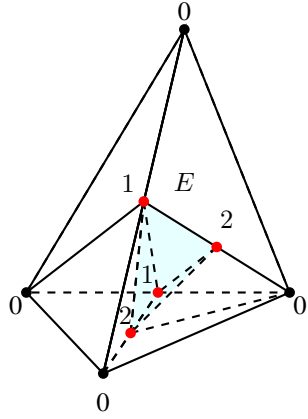
\begin{figure}[]
\centering
\begin{tikzpicture}[scale = 0.7]
\fill[bubbles] plot coordinates{(1,1.25,1) (-0.5, 0,0) (0, 2.5,2) (1,1.25,1) };
\fill[bubbles] plot coordinates{(-0.25,0,2) (-0.5, 0,0) (0, 2.5,2) (-0.25,0,2) };
\fill[bubbles] plot coordinates{(-0.25,0,2) (1,1.25,1) (0, 2.5,2) (-0.25,0,2) };
\fill[bubbles] plot coordinates{(-0.25,0,2) (1,1.25,1) (-0.5, 0,0) (-0.25,0,2) };
\draw [thick] plot coordinates{ (2,0,0)  (0,5,0) (0,0,4)  (2,0,0)} ;
\draw [thick] plot coordinates{ (-3, 0, 0) (0,5,0) (0,0,4) (-3, 0, 0)};
\draw [thick,dashed] plot coordinates{ (-3, 0, 0)  (0,5,0)  (2,0,0) (-3, 0, 0)} ;
\draw [thick] plot coordinates{ (0, 5, 0)  (0,0,4) } ;
\draw [thick] plot coordinates{(0, 2.5,2) (-3,0,0) } ;
\draw [thick] plot coordinates{(0, 2.5,2) (2,0,0) } ;
\draw [thick,dashed] plot coordinates{(0, 2.5,2)  (-0.5,0,0) } ;
\draw [thick,dashed] plot coordinates{(0, 0,4)  (-0.5,0,0) } ;
\draw [thick,dashed] plot coordinates{(0, 0,4)  (-0.5,0,0) } ;
\draw [thick,dashed] plot coordinates{(-0.25,0,2) (0, 2.5,2) } ;
\draw [thick,dashed] plot coordinates{(-0.25,0,2) (2, 0,0) } ;
\draw [thick,dashed] plot coordinates{(-0.25,0,2) (0, 2.5,2) } ;
\draw [thick, dashed] plot coordinates{ (1,1.25,1)  (-0.5, 0,0)} ;
\draw [thick, dashed] plot coordinates{ (1,1.25,1)  (-0.25,0,2) } ;
\draw [very thick, red]plot [mark=*, only marks]coordinates{(0, 2.5,2) (-0.5,0,0) (-0.25,0,2) (1,1.25,1)};
\draw [very thick]plot [mark=*, only marks]coordinates{(2, 0,0) (0,5,0) (0,0,4) (-3,0,0)};
\coordinate [label= 0] (z1) at (2.2,-0.6 ,0);
\coordinate [label= 0] (z1) at (-3.2,-0.6 ,0);
\coordinate [label= 0] (z1) at (0,5 ,0);
\coordinate [label= 0] (z1) at (0,-0.9,4);
\coordinate [label= $1$] (z1) at (-0.3,2.5,2);
\coordinate [label= $1$] (z1) at (-0.7,0,0);
\coordinate [label= $2$] (z1) at (-0.35,0,2) ;
\coordinate [label= $2$] (z1) at (1,1.25,1-0.5) ;
\coordinate [label= $E$] (z1) at (0,1.75,0) ;
\end{tikzpicture}
\caption{Black and red points represent the proper and the hanging nodes respectively. The light blue element has only hanging nodes as its vertices. Notice that this is obtained by a further refinement of the last element of Figure \ref{Figure:GlobalIndex}.}
\label{Figure:Example_all_hangings}
\end{figure}

The second preliminary property involves the interpolation error $v-\Ical_E  v$ on the boundary of an element $E\in \mesh$, for any function $v \in \Vmesh$. 
Let first $L$ be one of the six {\textit {straight sides}} of the element $E$.
If $L$ has not been refined, by definition $v|_L$ is a polynomial of degree one, and it is fully determined by the values of $v$ at the endpoints of $L$, which coincides with two vertices of $E$.
It results that $(v-\Ical_E v)|_L=0$. 
On the other hand, if $L$ contains one hanging node $\tilde{\x}$, $(\Ical_E v) (\tilde{\x}) = \frac{1}{2} v(\x')+ \frac{1}{2} v(\x'')$, with $\x'$ and $\x''$ the endpoints of $L$. This suggests the definition, firstly given in \cite{AVEMstabfree}, of the function of the \textit{hierarchical details} $d$ of $v$, as
    \begin{align}\label{def:difference_dx}
        d(v,\z)= \begin{cases}
            v(\z)&\text{if }\z \in \propers_E,\\
            v(\z) -\frac{1}{2}\left(v(\z')+ v(\z'')\right) &\text{if }\z\in \hangings_E,
        \end{cases}
    \end{align}
  where $\{\z',\z''\}$ is the set of endpoints of the edge having $\z$ as midpoint.  The links of this function with the global interpolation error has been clarified in \cite[Lemma 7.1 and Corollary 7.2]{AVEMstabfree} and we here present them in the three-dimensional case, providing the necessary changes in the proofs. 
  \begin{lemma}[global interpolation error vs hierarchical details]\label{lemma:global_vs_hierarchical_details}
      If $\mesh$ is a $\Lambda$-admissible partition, it holds
    \begin{equation*}
     \normHT{v-\Icalmesh v}^2 \simeq\,h\, \sum_{\x \in \hangings} d^2(v,\x), \quad \forall v \in \Vmesh,
    \end{equation*}
    where the hidden constants depend on $\Lambda$.
  \end{lemma}
\begin{proof}
Let $E\in \mesh$ be an element. From \eqref{e:stab_equivalence}, and the definition of $\Ical_E$, i.e. $\Ical_E v(\bm{x}) = v(\bm{x})$, $\forall \bm{x}\in\propers_E$,  we have the following equivalence of norms, $\forall v \in\V_E$ 
    \begin{align}\label{Hconstants1}
        \normHE{v - \Ical_E v}^2 \simeq s_E(v-\Ical_E v, v-\Ical_E v)= h_E\,\sum_{\x \in \hangings_E}\abs{v(\x)- \Ical_E v(\x)}^2.
    \end{align}
    We now define the subspace of $\V_E$, $\Tilde{\V}_E\coloneqq\left\{w \in \V_E: w(\x)=0,\;\forall \x \in \mathcal{P}_E\right\}$, observing that, in particular, $v-\Ical_E v \in \Tilde{\V}_E$. On this finite-dimensional space,  we can consider these two equivalent norms
    \begin{equation}\label{Hconstants2}
        h_E\sum_{\x \in \hangings_E} d^2(w,\x) \simeq h_E\sum_{\x \in \hangings_E} \abs{w(\x)}^2, \quad w \in \Tilde{\V}_E.
    \end{equation}
     We remark that the hidden constants depend on $\Lambda$, since, from  Assumption \ref{ass:LambdaMax}, both the dimension of the space $\Tilde{\V}_E$ and the number of the possible patterns of hanging nodes on $\partial E$ (from which the function $d(w,\bm{x})$ depends) can be bounded by $\Lambda$. Taking now $ (v - \Ical_E v) \in \Tilde{\V}_E$, and using \eqref{Hconstants1} and \eqref{Hconstants2} we have that
          \begin{align*}
          \normHE{v-\Ical_E v}^2 \simeq h_E\sum_{\x \in\hangings_E} d^2(v,\x),\quad\forall v\in\V_E.
      \end{align*}
      Finally, summing over all the elements of the partition, we conclude the proof.
\end{proof}
We now move to the space $\Vmeshz$, defined in \eqref{def:V_meshz}, where we set the Lagrange interpolation operator   
\begin{equation*}
    \Icalmeshz:\Vmesh \rightarrow \Vmeshz,
\end{equation*}
such that $(\Icalmeshz v)(\x) = v(\x)$, $\forall \x \in \propers$. 
We recall that a function $v\in \Vmeshz$ is a polynomial of degree one in three dimensions, thus it can be fully determined by the values at the four proper nodes $\propers_E$ at each element $E$. 
On the other hand, the value of $v$ at a hanging node $\x \in \hangings_E$ is determined as $(\Icalmeshz v)(\x) = \frac{1}{2} (\Icalmeshz v)(\x') + \frac{1}{2} (\Icalmeshz v)(\x'')$, where $\{\bm{x'}, \bm{x''}\}$ are the endpoints of the edge for which $\x$ is the midpoint, in a recursive way.

\begin{lemma}\label{lemma:Imesh_delta}
    Let $\mesh$ be a $\Lambda$-admissible partition and $\delta(v, \x) : = v(\x)- (\Icalmeshz v)(\x)$, it holds
    \begin{equation*}
    \normHT{\Icalmesh v - \Icalmeshz v}^2 \simeq h \sum_{\x\in\hangings}\delta^2(v,\x)\quad \forall v \in \Vmesh.
    \end{equation*}
\end{lemma}
\begin{proof}
    We fix an element $E\in\mesh$. If $v\in\Vmesh$, $\Ical_E v$ and $(\Icalmeshz v)|_E$ differs only at the vertices of $E$. Remarking that $\Ical_E v(\x)=v(\x)$ $\forall \x \in \propers_E$ we get
    \begin{equation*}
    \normHE{\Ical_E v -\Icalmeshz v}^2\simeq h_E \sum_{\x\in \propers_E} \abs{\left(\Ical_E v - \Icalmeshz v\right)(\x)}^2= h_E \sum_{\x\in \propers_E} \abs{\left(v - \Icalmeshz v\right)(\x)}^2.
    \end{equation*}
    We conclude the proof by summing up all the elements of the partition and recalling that $(\Icalmeshz v)(\x) = v(\x)$, $\forall \x \in \propers$, and that $\bigcup_{E\in\mesh}\propers_E = \nodes$, since any node is a vertex for at least one element, 
    \begin{equation*}
    \normHT{\Ical_E v -\Icalmeshz v}^2\simeq h \sum_{\x\in \nodes}  \abs{\left(v - \Icalmeshz v\right)(\x)}^2\simeq h \sum_{\x\in \hangings} \abs{\left(v - \Icalmeshz v\right)(\x)}^2.
    \end{equation*}
\end{proof}
\begin{prop}[comparison between interpolation operators] \label{prop:comparisonBetweenInterp}
Given a $\Lambda$-admissible partition $\mesh$ it holds
\begin{equation*}\label{Comparison}
    \normH{v -\Icalmeshz v}\lesssim \normHT{v -\Icalmesh v} \qquad \forall v \in \Vmesh,
\end{equation*}
where the hidden constant is independent of $\mesh$.
\end{prop}
\begin{proof}
    We first notice that using the triangle inequality it holds $\forall v \in \Vmesh$
    \begin{equation*}
        \normH{v- \Icalmeshz v} = \normHT{v -\Icalmeshz v}\le \normHT{v - \Icalmesh v}+\normHT{\Icalmesh v -\Icalmeshz v},
    \end{equation*}
    then it is enough to bound the last term, i.e.
    \begin{equation*}
        \normHT{\Icalmesh v -\Icalmeshz v}\lesssim \normH{v- \Icalmesh v} .
    \end{equation*}
    Thanks to Lemma \ref{lemma:global_vs_hierarchical_details} and Lemma \ref{lemma:Imesh_delta}, the proof reduces to show that
    \begin{equation*}
       h \sum_{\x\in\hangings} \delta^2(v, \bm{x}) \lesssim h \sum_{\x\in\hangings} d^2(v, \bm{x}), \quad \quad \quad \forall v \in \Vmesh,
    \end{equation*}
   or, in an equivalent form, simplifying the $h$, \begin{equation}\label{eq:delta_less_d}
          \norm{\bm{\delta}}_{l^2(\hangings)}\lesssim \norm{\bm{d}}_{l^2(\hangings)},
    \end{equation}
    where $\bm{\delta}=(\delta(\x))_{\x \in \hangings}\coloneqq (\delta (v, \x))_{\x \in \hangings}$ and $\bm{d} = (d(\x))\coloneqq(d(v,\x))_{\x \in \hangings}$.
    Given $\x \in \hangings$, $\{\x',\x'\}$ are the endpoints of the edge with $\x$ as midpoint, and exploiting \ref{def:difference_dx}, we have that
    \begin{align*}
        \delta(\x) &= v(\x) - \Icalmeshz v(\x) = v(\x)- \frac{1}{2}\Icalmeshz(\x') -\frac{1}{2}\Icalmeshz(\x'') \\&=
        v(\x) -\frac{1}{2}\left(v(\x') + v(\x'')\right) + \frac{1}{2}\left(v(\x') - \Icalmeshz v(\x')\right) + \frac{1}{2}\left(v(\x'') - \Icalmeshz v(\x'')\right) \\
        &=d(\x) +\frac{1}{2}\delta(\x') + \frac{1}{2}\delta(\x'').
    \end{align*}
    We can define a matrix $\bm{W}:l^2(\hangings)\rightarrow l^2(\hangings)$ such that $\bm{\delta} =  \bm{W}\bm{d}$. Then, \eqref{eq:delta_less_d} verifies if 
    \begin{equation*}
        ||\bm{W}||_2\lesssim 1.
    \end{equation*}
    We organize the global indices $\lambda\in[1,\Lambda_\mesh]$ and we define the set $\hangings_\lambda=\{\x\in \hangings:\lambda(\x)=\lambda\}$ and $\hangings =\bigcup_{1\le \lambda\le \Lambda_\mesh} \hangings_\lambda$ so that the  matrix $\bm{W}$ can be factorized in a block-wise manner as
    \begin{equation*}
    \bm{W} = \bm{W}_{\Lambda_\mesh} \bm{W}_{\Lambda_\mesh -1} \cdots\bm{W}_{1},
    \end{equation*}
    where $\bm{W}_1$ is the identity matrix since $\delta(\x') = \delta(\x'')=0$ and each matrix $\bm{W}_\lambda$ is obtained changing from the identity matrix the rows of block $\lambda$, i.e.,
    %
\begin{linenomath*}
\begin{align*}
\bm{W}_\lambda=
\begin{matrix}
\begin{bmatrix} 1 & 0 &\cdots & \cdots& \cdots & \cdots& 0 \\
0 & 1 &0 & \cdots & \cdots & \cdots& 0 \\
\vdots & \vdots & \ddots  &{}&{}  &{}& \vdots \\
\vdots & \vdots &{} & \ddots  &{}&{}  & \vdots \\
\hline
\vdots &{1/2}  &{} &{1/2} &\vline \; \;1\;\;\vline & {} & \vdots \\
\hline
\vdots &\vdots &{}&{}  &{}& \ddots  &{}\vdots \\
0& 0& \cdots  &{} &{} &{}&1
\end{bmatrix} &\begin{matrix}
\vspace{0.5cm}
    {}\\
    {}\\
    {}\\
    {}\\
    \leftarrow\; {\text {block of hanging nodes of index }}\; \lambda\\
    {}\\
    {}
\end{matrix}
\end{matrix}.
\end{align*}
\end{linenomath*}
    
    We can now apply the H\"older inequality as $\|\bm{W}_\lambda\|_2\le \left(\|\bm{W}_\lambda\|_1 \|\bm{W}_\lambda\|_\infty\right)^{1/2}$,\\ $\|\bm{W}_\lambda\|_\infty~\le~ \frac{1}{2}~+\frac{1}{2}~+~1~ =2$, and $\|\bm{W}_\lambda\|_1\le \frac{1}{2} c +1$. In particular, for the second term we exploit the fact that the number of hanging nodes (and so the constant $c$) depends on the quality of the mesh, i.e. on the minimum three-dimensional angle required between each face of the partition.
    Thus, 
    \begin{equation*}
    \|\bm{W}\|_2\le \prod_{2\le \lambda\le \Lambda_\mesh}\|\bm{W}_\lambda\|_2 \le (c+2)^{\frac{\Lambda-1}{2}},
    \end{equation*}
    which ends the proof.
\end{proof}
We finally introduce  $\Icalmeshzt: \Vmesh \rightarrow \Vmeshz$ the classical Clément interpolation operator on $\Vmeshz$, which takes at each internal node the average target function on the support of the associated basis function. The following Lemma has been introduced in \cite[Lemma 3.3]{AVEMstabfree}, it is based on Proposition \ref{ScaledPoincare}, and it does not depend on the geometric dimension.
\begin{lemma}[Clément interpolation estimate]\label{lemma:ClementInterpolationEstimate} It holds
    \begin{equation*}
         \sum_{E \in \mesh}h^{-2}_E \normLE{v-\Icalmeshzt v}^2 \lesssim \normH{v}^2 \qquad \forall v \in \V.
    \end{equation*}
\end{lemma}
\begin{proof}
Let $\Icalmesht: \V \rightarrow \Vmesh$ be the Clément operator on the virtual-element space, see \cite{Clement} for the details. Given $v\in \V$, we denote by $\tilde{v}  =\Icalmesht v \in \Vmesh$. We write 
$ v - \Icalmeshzt v = (v - \vt) + (\vt- \Icalmeshzt \vt) +(\Icalmeshzt \vt - \Icalmeshzt v)$ and from the local stability of $\Icalmeshzt$ in $L^2$, we get
\begin{equation*}
\begin{split}
    \sum_{E\in \mesh}h^{-2}_E\normLE{v- \Icalmeshzt v}^2 \lesssim \sum_{E\in \mesh}h^{-2}_E \normLE{v -\vt}^2 + h_E^{-2}\normLE{\vt -\Icalmeshzt \vt}^2 \\\lesssim \normH{v}^2 + \sum_{
E\in \mesh}h_E^{-2}\normLE{\vt - \Icalmeshzt \vt}^2.
\end{split}
\end{equation*}
Furthermore, from the invariance of $\Icalmeshzt$ in $\Vmeshz$, i.e. $\Icalmeshzt (\Icalmeshz \vt) = \Icalmeshz \vt$, and we have $\vt - \Icalmeshzt \vt = (\vt - \Icalmeshz \vt) + \Icalmeshzt(\Icalmeshz \vt - \vt)$. We can then conclude the proof using Proposition \ref{ScaledPoincare} and the stability of $\Icalmeshzt$ in $L^2$. 
\end{proof}

\section{A posteriori error analysis}\label{Sec:AposteroriErrorAnalysis}
In this section we prove the bound of the stabilization term by the residual.
We essentially rely on the validity of the scaled Poincaré inequality \ref{ScaledPoincare}, Proposition \ref{Comparison} and Lemma \ref{lemma:ClementInterpolationEstimate}.
We highlight that this analysis follows the same steps of the one presented in \cite{AVEMstabfree}, as it is independent on the geometric dimension of the problem.

Following the a posteriori error analysis for the VEM in \cite{Cangiani}, we introduce the internal residual associated with element $E\in\mesh$ as
\begin{equation*}
    \res(E; v, \data)\coloneqq f_E - c_E \Pi^\nabla_E v, \quad \forall\, v \in \V_E
\end{equation*}
and the jump residual over a face $F\coloneqq E_1 \cap E_2 \subset \Omega$ of two elements $E_1$, $E_2$
\begin{equation*}
    \jumpT(F; v, \data)\coloneqq [[ \K  \nabla \Pi^\nabla_\mesh  v]]_F = \left(\K_{E_1} \nabla \Pi^\nabla_{E_1} v|_{E_1}\right)\cdot \bm{n_1} + \left(\K_{E_2} \nabla \Pi^\nabla_{E_2} v|_{E_2}\right)\cdot \bm{n_2},
\end{equation*}
where $\bm{n_i}$ denotes the unit normal vector to $F$ that points outward to the element $E_i$, $i=1,2$. If the face $F\subset \partial \Omega$, we set $\jumpT(F; v, \data)\coloneqq 0$. We define the local estimator as
\begin{equation}\label{eq:localEstimator}
    \eta_\mesh^2(E;v,\data) \coloneqq h_E^2\normLE{\res(E; v,\data)}^2 + \frac{1}{2}\sum_{F\in \mathcal{F}_E}h_E  \norm{\jumpT(F; v, \data)}^2_{0, F},
\end{equation}
and the global estimator as the sum over all the elements of the partition
\begin{equation*}
\eta^2_\mesh (v, \data)\coloneqq \sum_{E\in \mesh}\eta_\mesh^2(E;v,\data), \quad \forall v \in \Vmesh.
\end{equation*}
With this error estimator, proceeding as in \cite{Cangiani},  and \cite[Proposition 4.1]{AVEMstabfree} the following upper and lower bounds are stated.
We remark that in the following propositions, the quality of $\mesh$ affects the constants. However, we have chosen to omit this dependence in the statements, as it is not the primary focus of this paper.
\begin{prop}[upper bound]\label{prop:upperBound} There exists a  $C_{apost}>0$ depending only on $\Lambda$ and $\data$, such that
\begin{equation*}
\normH{u- u_\mesh }^2 \le C_{apost} \left(\eta^2_\mesh (u_\mesh,\data) +S_\mesh(u_\mesh, u_\mesh)\right).
\end{equation*}
\end{prop}
\begin{prop}[lower bound]\label{prop:lowerBound} There exists a constant $c_{apost}>0$ depending only on $\Lambda$, such that
\begin{equation*}
c_{apost} \;\eta^2_\mesh (u_\mesh,\data) \le \normH{u- u_\mesh }^2 +S_\mesh(u_\mesh, u_\mesh) .
\end{equation*}
\end{prop}
The preliminary properties obtained in the previous section, allow us to obtain the bound of the stabilization term, operating exactly as  \cite[Proposition 4.4]{AVEMstabfree}. For the sake of completeness, we report the full proof, highlighting in particular the use of the scaled Poincaré inequality (Proposition \ref{ScaledPoincare}).
\begin{prop}[bound of the stabilization term by the residual]\label{Bound} There exists a constant $C_B$ depending only on $\Lambda$, such that
\begin{equation}\label{eq:boundS}
    \gamma^2 S_\mathcal{T}(u_\mathcal{T},u_\mathcal{T})\le C_B \eta^2_\mathcal{T}(u_\mathcal{T}, \mathcal{D}).
\end{equation}
\end{prop}
\begin{proof}
Let  $w\in\V^0_\mesh$ and $e_\mesh : =  u_\mesh -w$. From the definition of the bilinear form $\B_\mesh (\cdot, \cdot)$ \eqref{def:definition_BT}, we get
    \begin{align*}
        \gamma S_\mesh(u_\mesh,u_\mesh) &= \gamma S_\mesh(u_\mesh,e_\mesh ) \\&= \B_\mesh(u_\mesh,e_\mesh ) -a_\mesh\left(u_\mesh, e_\mesh \right)- m_\mesh \left(u_\mesh, e_\mesh \right)\\
        &= \Fcalmesh(e_\mesh ) -a_\mesh\left(u_\mesh, e_\mesh \right)- m_\mesh \left(u_\mesh, e_\mesh \right)\\
        &= \sum_{E\in\mesh} \int_E(f_E- c_E\Pi^\nabla_E u_\mesh)\Pi^\nabla_E e_\mesh -\sum_{E\in\mesh}\int_E \K_E \left(\nabla \Pi^\nabla_E u_\mesh \right)\cdot\left(\nabla \Pi^\nabla_E e_\mesh \right).
    \end{align*}
We apply \eqref{cond_PiNabla1} and the integration by part for the second term, noting that $\K_E$ is constant on $E$ and $\Delta \Pi^\nabla_E u_\mesh = 0$, which yields to
\begin{align*}
    \gamma S_\mesh(u_\mesh,u_\mesh) &= \sum_{E\in\mesh} \int_E(f_E- c_E\Pi^\nabla_E u_\mesh)\Pi^\nabla_E e_\mesh -\sum_{E\in\mesh}\int_{\partial E} \bm{n} \cdot\left(\K_E \nabla \Pi^\nabla_E u_\mesh \right) e_\mesh.
\end{align*}
Employing now the continuity of $\Pi^\nabla_E$ as stated in \eqref{eq:continuity_PiNabla}, we get
 \begin{align*}   
    \gamma S_\mesh(u_\mesh,u_\mesh)\le  &\sum_{E\in\mesh}h_E\normLE{ \res(E;  u_\mesh, \data)}h_E^{-1}\left(\normLE{e_\mesh}+ h_E\normHE{e_\mesh}\right)\\
    &+\frac{1}{2}\sum_{E\in \mesh}\sum_{F\in \mathcal{F}_E}h_E^{1/2}\norm{j_\mesh(F; u_\mesh, \data)}_{0,F} h_E^{-1/2}\norm{e_\mesh}_{0,F}.
\end{align*}
    For any $\delta >0$, the latter inequality brings to
    \begin{equation*}
    \gamma S_\mesh(u_\mesh, u_\mesh) \le \frac{1}{2 \delta} \eta^2_\mesh(u_\mesh, \data) + \frac{\delta}{2} \Phi_\mesh(e_\mesh)
    \end{equation*}
    where
    \begin{align*}
            \Phi_\mesh(e_\mesh)&=\sum_{E\in \mesh}\left(h^{-2}_E \normLE{e_\mesh}^2 + \normHE{e_\mesh}^2 + \sum_{F\in \mathcal{F}_E}h_E^{-1}\norm{e_\mesh}^2_{0,F}\right)\\ &\lesssim \sum_{E\in\mesh}\left(h_E^{-2}\normLE{u_\mesh -w}^2 +\normHE{u_\mesh - w}^2\right).
    \end{align*}
Choosing now $w = \Icalmeshz u_\mesh$, we can apply the scaled Poincaré~\eqref{ScaledPoincare}, obtaining 
\begin{equation*}
    \Phi_\mesh(u_\mesh - \Icalmeshz u_\mesh)\lesssim \abs{u_\mesh - \Icalmeshz u_\mesh}^2_{1,\Omega}. 
\end{equation*}
Employing Proposition \ref{prop:comparisonBetweenInterp} and \eqref{eq:S_mesh_Hnorm}, we get $\Phi_\mesh(u_\mesh - \Icalmeshz u_\mesh )\lesssim S_\mesh(u_\mesh, u_\mesh)$, which concludes the proof for a suitable choice of parameter $\delta$.
\end{proof}
The propositions~\ref{prop:upperBound}, \ref{prop:lowerBound} and \ref{Bound}, bring to the stabilization-free upper and lower bounds.
\begin{cor}[stabilization-free a posteriori error estimates] There exist two constants $C_L>~0$ and $C_U>0$ such that
\begin{equation}\label{eq:cor_boundeta}
C_L\eta^2_\mathcal{T}(u_\mathcal{T}, \mathcal{D})\le \normH{u- u_\mesh }^2 \le C_U \eta^2_\mathcal{T}(u_\mathcal{T}, \data).
\end{equation}   
\end{cor}
\section{The module {\sf GALERKIN}}\label{Sec:TheModuleGalerkin}
In analogy of \cite{AVEMstabfree}, we present the extended adaptive algorithm, {\sf AVEM},  consisting of an outer loop  made of two modules {\sf DATA} and {\sf GALERKIN}.
The module {\sf DATA} approximates at each step the original smooth $\data$ of the variational problem \eqref{Variational_Problem} with piecewise constants on the elements of the partition with prescribed accuracy. 
The {\sf GALERKIN} module takes as input a $\Lambda$-admissible partition, the approximated data and a given tolerance $\epsilon$, it iterates on the loop  
\begin{equation}\label{eq:Interations}
    {\sf SOLVE} \rightarrow {\sf ESTIMATE} \rightarrow {\sf MARK} \rightarrow {\sf REFINE},
\end{equation}
producing a sequence of $\Lambda$-admissible  partitions and the associated Galerkin approximate solution $u_k$, and stopping when the condition $\eta(u_k, \data)<\epsilon$ is satisfied. 
We here discuss only the {\sf GALERKIN} module, leaving the complexity and the quasi-optimality of the AVEM algorithm in a forthcoming paper, as done in \cite{AVEMConvergenceOptimality} for dimension two.

In particular, we highlight that in our analysis the modules in loop \eqref{eq:Interations} involve:
\begin{itemize}
    \item the D\"{o}rfler criterion \cite{Dorfler} in the {\sf MARK} step, i.e., given an input a parameter $\theta\in(0,1)$, it produces an almost minimal set $\mathcal{M}\subset \mesh$ such that
    \begin{equation*}
        \theta\; \eta^2_\mesh(u_\mesh, \data)\le \sum_{E\in\mathcal{M}}\eta^2_\mesh(E;u_\mesh, \data);
    \end{equation*}
    \item the newest-vertex bisection (NVB) refinement strategy for tetrahedral-shape elements in the  {\sf REFINE} step \cite{NoSiVe:09,StevensonNVB}. 
\end{itemize}
We also remark that in the {\sf REFINE} step a further module, {\sf MAKE\_ADMISSIBLE}, is needed to reestablish the $\Lambda$-admissibility.

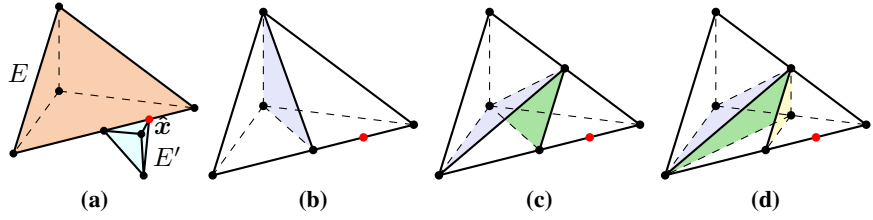
\begin{figure}[]
\centering
 \begin{subfloat}[]{
 \label{fig:MAKE_ADMISSIBLE_two_elements}
\begin{tikzpicture}[scale = 0.45]
\coordinate (A) at (0,0.5 ,0);
\coordinate (B) at (4,0 ,0);
\coordinate (C) at (0,3,0);
\coordinate (C1) at (2.5,-2,0);
\coordinate (D) at (0,0,3.5);
\coordinate (M1) at (2,0,1.75);
\coordinate (M2) at (3,0,0.875);
\coordinate (M3) at (3.2,0,2);

\fill[bubbles] plot coordinates{ (M2) (C1)(M1)(M2)} ;
\fill[apricot] plot coordinates{  (C)  (B) (D) } ;
\draw [thin,dashed] plot coordinates{  (A) (B)} ;
\draw [thick] plot coordinates{ (C) (D) };
\draw [thick] plot coordinates{ (B) (C)} ;
\draw [thick] plot coordinates{ (D)(B)  } ;
\draw [thin,dashed] plot coordinates{ (A) (C) } ;
\draw [thin,dashed] plot coordinates{  (A) (D) } ;

\draw [thick] plot coordinates{ (M2) (C1)(M1)(M3) (M2)} ;
\draw [thick] plot coordinates{ (C1)(M3)} ;

\draw [very thick]plot [mark=*, only marks]coordinates{(A) (B)(C)(D)  (M1) (M3) (C1) };
\draw [very thick,red]plot [mark=*, only marks]coordinates{(M2)  };
\coordinate [label= $\bm{\hat{x}}$] (z1) at (3.4,-0.7,0.875) ;
\coordinate [label= $E'$] (z1) at (3.2,-2 ,0);
\coordinate [label= $E$] (z1) at (-1.2,0.5 ,0);
\end{tikzpicture}}
\end{subfloat}
\begin{subfloat}[]{
\label{fig:MAKE_ADMISSIBLE_one_refinement}
\begin{tikzpicture}[scale = 0.5]
\coordinate (A) at (0,0.5 ,0);
\coordinate (B) at (4,0 ,0);
\coordinate (C) at (0,3,0);
\coordinate (D) at (0,0,3.5);
\coordinate (M1) at (2,0,1.75);
\coordinate (M2) at (3,0,0.875);
\coordinate (M3) at (3.2,0,2);
\fill[lavender(web)] plot coordinates{ (M1) (C)  (A) } ;
\draw [thin,dashed] plot coordinates{  (A) (B)} ;
\draw [thick] plot coordinates{ (C) (D) };
\draw [thick] plot coordinates{ (B) (C)} ;
\draw [thick] plot coordinates{ (D)(B) } ;
\draw [thin,dashed] plot coordinates{ (A) (C) } ;
\draw [thin,dashed] plot coordinates{  (A) (D) } ;


\draw [thick] plot coordinates{ (C)(M1)} ;
\draw[thin, dashed] plot coordinates{(A) (M1)};
\draw [very thick]plot [mark=*, only marks]coordinates{(A) (B)(C)(D)  (M1) };
\draw [very thick,red]plot [mark=*, only marks]coordinates{(M2)  };
\end{tikzpicture}
}\end{subfloat}
\begin{subfloat}[]{
\label{fig:MAKE_ADMISSIBLE_two_refinement}
\begin{tikzpicture}[scale = 0.5]
\coordinate (A) at (0,0.5 ,0);
\coordinate (B) at (4,0 ,0);
\coordinate (C) at (0,3,0);
\coordinate (D) at (0,0,3.5);
\coordinate (F) at (2,1.5,0);
\coordinate (M1) at (2,0,1.75);
\coordinate (M2) at (3,0,0.875);
\coordinate (M3) at (3.2,0,2);
\fill[grannysmithapple] plot coordinates{ (F)  (A)(M1) } ;
\fill[lavender(web)] plot coordinates{ (A) (F)  (D) } ;
\draw [thin,dashed] plot coordinates{  (A) (B)} ;
\draw [thick] plot coordinates{ (C) (D) };
\draw [thick] plot coordinates{ (B) (C)} ;
\draw [thick] plot coordinates{ (D) (B)} ;
\draw [thin,dashed] plot coordinates{ (A) (C) } ;
\draw [thin,dashed] plot coordinates{  (A) (D) } ;


\draw [thin, dashed] plot coordinates{ (A)(F)} ;
\draw [thin, dashed] plot coordinates{ (M1)(A)} ;
\draw [thick] plot coordinates{ (F)(D)} ;
\draw [thick] plot coordinates{ (F)(M1)} ;
\draw [very thick]plot [mark=*, only marks]coordinates{(A) (B)(C)(D)  (M1)  (F) };
\draw [very thick,red]plot [mark=*, only marks]coordinates{(M2)  };
\end{tikzpicture}
}\end{subfloat}
\begin{subfloat}[]{
\label{fig:MAKE_ADMISSIBLE_three_refinement}
\begin{tikzpicture}[scale = 0.5]
\coordinate (A) at (0,0.5 ,0);
\coordinate (B) at (4,0 ,0);
\coordinate (C) at (0,3,0);
\coordinate (D) at (0,0,3.5);
\coordinate (F) at (2,1.5,0);
\coordinate (G) at (2,0.25,0);
\coordinate (M1) at (2,0,1.75);
\coordinate (M2) at (3,0,0.875);
\coordinate (M3) at (3.2,0,2);
\fill[grannysmithapple] plot coordinates{ (F)  (G)(D) };
\fill[lemonchiffon] plot coordinates{ (F) (G)  (M1) } ;
\fill[lavender(web)] plot coordinates{ (F)  (A)(D) } ;
\draw [thin,dashed] plot coordinates{  (A) (B)} ;
\draw [thick] plot coordinates{ (C) (D) };
\draw [thick] plot coordinates{ (B) (C)} ;
\draw [thick] plot coordinates{ (D) (B) } ;
\draw [thin, dashed] plot coordinates{ (M1) (M2)} ;
\draw [thin,dashed] plot coordinates{ (A) (C) } ;
\draw [thin,dashed] plot coordinates{  (A) (D) } ;


\draw [thin, dashed] plot coordinates{ (A)(F)} ;
\draw [thin, dashed] plot coordinates{ (G)(F)} ;
\draw [thin, dashed] plot coordinates{ (G)(D)} ;
\draw [thin, dashed] plot coordinates{ (G)(M1)} ;
\draw [thick] plot coordinates{ (F)(M1)} ;
\draw [thick] plot coordinates{ (F)(D)} ;
\draw [very thick]plot [mark=*, only marks]coordinates{(A) (B)(C)(D)  (M1) (F) (G) };
\draw [very thick,red]plot [mark=*, only marks]coordinates{(M2)  };
\end{tikzpicture}
}\end{subfloat}
\caption{
Two elements $E$ (peach) and $E'$ (light blue) sharing an edge and the node $\bm{\hat{x}}$ (red) that violates Assumption \ref{ass:LambdaMax} (a). If $\bm{\hat{x}}$ does belong to the newest-edge, one refinement of $E$ restores the $\Lambda$-admissibility (b).  If $\bm{\hat{x}}$ does not belong to the newest-edge, two and three refinements of $E$ are needed (c) and (d), respectively. Purple, green and yellow are the new faces separating the tetrahedral-shape elements. In (b), (c), and (d) $E'$ has not been represented for graphic purposes.}
\label{Figure:MAKE_ADMISSIBLE}
\end{figure}

Indeed, let us suppose that after a refinement a node $\bm{\hat{x}}$ of an element $E'$ has global $\lambda(\bm{\hat{x}}) = \Lambda +1$, violating the Assumption \ref{ass:LambdaMax}. By definition of global index of a node, $\bm{\hat{x}}$ is a hanging node for at least an other tetrahedral-shape element $E$ (see Figure \ref{fig:MAKE_ADMISSIBLE_two_elements}). Two possible strategies can be adopted to restore the $\Lambda$-admissibility, which depend on the position of the node $\bm{\hat{x}}$. If $\bm{\hat{x}}$ belongs to the newest-edge of the element $E$ (Figure \ref{fig:MAKE_ADMISSIBLE_one_refinement}), one refinement of element $E$ brings to the condition $\lambda(\bm{\hat{x}})\le \Lambda$. 
On the other hand, if the edge 
 in which $\bm{\hat{x}}$ lays cannot be immediately refined (as in Figure \ref{fig:MAKE_ADMISSIBLE_two_refinement} and \ref{fig:MAKE_ADMISSIBLE_three_refinement}), then the NVB produces, respectively, three and four new elements.

\subsection{The convergence of {\sf GALERKIN}}

We discuss the reduction of the a posteriori error estimator under the refinement, to prove the convergence of {\sf GALERKIN} from a mesh $\mesh$ to the refined mesh $\mesh_*$. 
We consider an element $E \in \mesh$ that splits into two elements $E_i \in \mesh_*$, $i=1,2$, with the introduction of a new face $F=E_1\cap E_2$.
The value of a function $v\in \Vmesh$ is known at all the nodes and edges of $E$, and then at the newest-vertex produced by the refinement. The new edges introduced do not contain any internal nodes, and then also the value of $v$ at $\partial F$ is known. We can then deduce that  $v$ is known at all the nodes on $\partial E_1$ and $\partial E_2$. 
We then associate to any $v \in \Vmesh$ a function $v_*\in \V_{\mesh_*}$  such that $v|_F = v_*|_F$. 

We denote by $\eta_{\mesh_*}(E; v_*, \data)$ the sum of the local estimators on the two new elements, i.e.
\begin{equation*}
\begin{split}
    \eta^2_{\mesh_*}(E; v_*, \data)\coloneqq \sum^2_{i=1}\eta^2_{\mesh_*}(E_i;v_*, \data)\quad \quad\quad \quad\quad \quad\quad \quad\quad\quad \quad\quad \quad\quad\quad \quad\quad \quad \quad \quad\\  =\sum^2_{i=1}h^2_{E_i}\norm{r_{\mesh_*}(E; v_*, \data)}^2_{0,E_i} + \sum^2_{i=1}\frac{h_{E_i}}{2}\sum_{F\in \partial \mathcal{F}_{E_i}} \norm{j_{\mesh_*}(F; v_*, \data)}^2_{0, F},
\end{split}
\end{equation*}
where $h_{E_i}=\frac{1}{\sqrt[3]{2}}h_E$, since the volume of $E$ halves after the refinement. 
Furthermore, we remark that $\K_E|_{E_i} = \K_{E_i}$, $c_E|_{E_i}= c_{E_i}$, and $f_{E}|_{E_i}= f_{E_i}$. 
As a first step, we analyze the behavior of the estimator under the refinement in analogy to \cite[Lemma 5.2] {AVEMstabfree}.

\begin{lemma}[local residual under refinement]
    Let $E$ be an element of the partition $\mesh$ which is split into two elements $E_1$ and $E_2$ in $\mesh_*$. There exists a strictly positive constant $\mu<1$  such that $\forall v \in \Vmesh$
    \begin{equation*}
    \eta_{\mesh_*}(E;v_*,\data)\le \mu\, \eta_\mesh(E;v, \data)+ c_{er,1} S^{1/2}_{\mesh(E)}(v,v),
    \end{equation*}
    where $c_{er,1}>0$ and  
    \begin{equation*}
        \mesh(E)\coloneqq \{E'\in \mesh: \mathcal{F}_E \cap \mathcal{F}_{E'}\neq \emptyset\}, \ S_{\mesh(E)}(v,v)\coloneqq\sum_{E'\in \mesh(E)}S_{E'}(v,v).
    \end{equation*}
\end{lemma}
\begin{proof}
    As a first step, we write the residual on $E$ and on its children $E_i$, $i=1,2$,
    \begin{equation*}
        r_E = f-c\,\Pi^\nabla_E v, \quad \quad \quad r_{E_i} =  f-c\,\Pi^\nabla_{E_i} v_*,
    \end{equation*}
    where we have dropped the $E$ from the piecewise constant data $f$ and $c$.
    Note that $r_{E_i}= r_E +c \left(\Pi^\nabla_E v -\Pi^\nabla_{E_i} v_*\right)$ and, given $\epsilon>0$
    \begin{equation*}
        \sum^2_{i=1}h^2_{E_i}\norm{r_{E_i}}_{0,E_i}^2\le \sum_{i=1}^2 h^2_{E_i}\left((1+ \epsilon)\norm{r_E}^2_{0,E_i} +\left(1 + \frac{1}{\epsilon}\right)\norm{c(\Pi^\nabla_E v - \Pi^\nabla_{E_i} v_*)}^2_{0,E_i}\right).
    \end{equation*}
    Recalling that $h_{E_i}=\frac{1}{\sqrt[3]{2}}h_E$, we get
    \begin{equation}\label{eq:lemma_residual}
        \sum^2_{i=1}h^2_{E_i}\norm{r_{E_i}}_{0,E_i}^2\le
        \frac{1+\epsilon}{\sqrt[3]{4}}h^2_E \norm{r_E}^2_{0,E}+ h^2_E \left(1 + \frac{1}{\epsilon}\right)\frac{c^2}{\sqrt[3]{4}}\sum^2_{i=1}\norm{\Pi^\nabla_E v - \Pi^\nabla_{E_i} v_*}^2_{0,E_i}
    \end{equation}
From the triangle inequality, one has
\begin{equation*}
    \norm{\Pi^\nabla_E v - \Pi^\nabla_{E_i} v_*}^2_{0,E_i}\le  \norm{v-\Pi^\nabla_E v}^2_{0, E} + \sum_{i=1}^2 \norm{v_*- \Pi^\nabla_{E_i}v_*}^2_{0,E_i},
\end{equation*}
and applying the minimality property of $\Pi^\nabla$ \eqref{diseq:PiNabla}, we get 
\begin{equation*}
    \norm{v-\Pi^\nabla_E v}^2_{0,E} \lesssim \normHE{v- \Pi^\nabla_E v}^2\le \normHE{v-\Ical_E v}^2,
    \end{equation*}
    \begin{equation*}    
    \sum^2_{i=1}\norm{v_* -\Pi^\nabla_{E_i} v_*}^2_{0,E_i}\lesssim \sum^2_{i=1}\abs{v_*-\Pi^\nabla_{E_i} v_*}^2_{1,E_i}\le\normHE{v-\Pi^\nabla_E v}^2 \le \normHE{v-\Ical_E v}^2.
\end{equation*}
If we set, for instance, $\epsilon = \frac{1}{2}$, and $\mu\coloneqq \frac{1+\epsilon}{\sqrt[3]{4}}=\frac{3}{2\sqrt[3]{4}}<1$, the inequality~\eqref{eq:lemma_residual} becomes
\begin{equation*}
    \sum^2_{i=1}h_{E_i}^2\norm{r_{E_i}}^2_{0,E_i}\le \mu\, h^2_E\normLE{r_E}^2 + C h_E^2 \normHE{v- \Ical_E}^2,
\end{equation*}
with a suitable constant $C>0$. Employing \eqref{eq:S_mesh_Hnorm} we get
\begin{equation*}
    \sum^2_{i=1}h_{E_i}^2\norm{r_{E_i}}^2_{0,E_i}\le \mu\, h^2_E\normLE{r_E}^2 + C h_E^2 S_E(v,v).
\end{equation*}
For the jump term, we have that for any $F\in \partial E$,  
\begin{equation*}
    j_{\mesh_*}(F;v_*) = \jumpT(F; v) +\left(j_{\mesh_*}\left(F; v_* \right) - j_{\mesh} (F; v)\right).
\end{equation*}
Fix $\epsilon>0$, we operate as we have done for the residual term, i.e.,
\begin{equation*}
    \sum^2_{i=1} \sum_{F\in  \mathcal{F}_{E_i} }h_{E_i}\norm{j_{\mesh_*}(F; v_*)}^2_{0, F}\le (1+\epsilon)\, I + \left(1 + \frac{1}{\epsilon}\right)\,II
\end{equation*}
where 
\begin{equation*} 
    I =  \sum^2_{i=1} \sum_{F\in \mathcal{F}_{E_i} }h_{E_i} \norm{\jumpT(F;v)}^2_{0,F},\ 
    II = \sum^2_{i=1} \sum_{F\in \mathcal{F}_E}h_{E_i}\norm{j_{\mesh_*}(F;v_*)-j_\mesh (F;v)}^2_{0,F}.
\end{equation*} 
On the new face that is generated in the refinement the jump term is zero, therefore 
\begin{equation*}
    I\le \frac{1}{\sqrt[3]{2}}\sum_{F\in \mathcal{F}_E}h_E \norm{\jumpT(F;v)}^2_{0,F}.
\end{equation*}
For the second term, we start by denoting by $\mesh_*(E_i)\coloneqq\{E' \in \mesh_*:\mathcal{F}_{E_i}\cap \mathcal{F}_{E'}\neq \emptyset\}$ and the element $E_{i,F}\in \mesh_*(E_i)$ such that, for any face $F\in \mathcal{F}_{E_i}$, $F =  \partial E_i \cap \partial E_{i,F}$. Indicating by $\hat{E}_{i,F}$ the parent of $E_{i,F}$, we have
\begin{equation*}
\begin{split}
    \norm{j_{\mesh_*}(F;v_*) -j_\mesh (F;v)}_{0,F} =\norm{[[ \K \nabla \left(\Pi^\nabla_{\mesh_*} - \Pi^\nabla_{\mesh}\right) v]] }_{0,F}\quad \quad \quad \quad \quad \quad \quad \quad \quad\quad \quad \quad \quad \\\le \norm{\K_E \nabla \left(\Pi^\nabla_{E_i} -\Pi^\nabla_E\right)v}_{0,F} + \norm{\K_{\hat{E}_{i,F}} \nabla (\Pi^\nabla_{E_i, F} - \Pi^\nabla_{\hat{E}_i,F})v }_{0,F},
    \end{split}
\end{equation*}
Employing the trace inequality, we get
\begin{align*}
    II &\lesssim \sum^2_{i=1}\sum_{E'\in \mesh_*(E_i)} \norm{\nabla(\Pi^\nabla_{E'} -\Pi^\nabla_{\hat{E}'} ) v}^2_{0,E'}\\
    &\lesssim \sum^2_{i=1}\sum_{E'\in\mesh_*(E_i)}\left(\norm{\nabla v - \nabla \Pi_E^\nabla v}^2_{0,E'} + \norm{\nabla v - \nabla \Pi^\nabla_{\hat{E}'}v}^2_{0,E'}\right)
\end{align*}
From the minimality of the operators $\Pi^\nabla_{E'}$ and $\Pi^\nabla_{\hat{E}'}$, we get
\begin{equation*}
    II \lesssim \sum_{E'\in \mesh(E)}\norm{\nabla(v-\Ical_{E'}v)}^2_{0,E'}\lesssim \sum_{E'\in \mesh(E)} S_{E'} (v,v),
\end{equation*}
which brings to the conclusion of the proof, with a sufficiently small $\epsilon$.
\end{proof}
This result, in addition to the Lipschitz continuity of the estimator with respect to the function $v$ (whose proof can be found in \cite[Lemma 5.3]{AVEMstabfree} and does not depend on the dimension), leads to the following Proposition. 
\begin{prop}[estimator reduction]
    There exist three positive constants independent of $\mesh$, namely $\rho<1$, $C_{er,1}$, and $C_{er,2}$, such that it holds $\forall w \in \V_{\mesh_*}$
    \begin{equation}\label{eq:reduction}
    \eta^2_{\mesh_*}(\mesh_*;w,\data) \le \rho \eta^2_\mesh(\mesh; u_\mesh,\data) + C_{er,1} S_\mesh (u_\mesh,u_\mesh) + C_{er,2} \abs{u_\mesh -w}^2_{1,\Omega},
    \end{equation}
    where $u_\mesh \in \Vmesh$ is the Galerkin solution and $\mesh_*$ the refined partition obtained by applying the step {\sf REFINE}.
\end{prop}
The second fundamental Proposition that involves the bound of the stabilization term is the quasi-orthogonality property \cite[Corollary 5.8]{AVEMstabfree}. 

\begin{prop}[quasi-orthogonality of energy errors without stabilization]
Given $\delta\in \left(0, \frac{1}{4}\right)$, it holds the following
\begin{equation}\label{eq:quasi_orth}
\EnergyNorm{u -u_{\mesh_*}}^2 \le (1+ 4 \delta) \EnergyNorm{u - u_\mesh}^2 - \EnergyNorm{u_{\mesh_*}- u_\mesh}^2,
\end{equation}
where $u_{\mesh_*}\in \V_{\mesh_*}$ is the solution of problem \eqref{Discrete_Problem} on the refine mesh $\mesh_*$.
\end{prop}

These two properties
guarantee the convergence of the {\sf GALERKIN} module with the same steps of \cite[Theorem 5.1]{AVEMstabfree}.
\begin{thm}[contraction property of {\sf GALERKIN}]
    Let $u_{\mesh}\in \Vmesh$ the Galerkin solution, and $\mathcal{M}\subset \mesh$ the set of marked elements obtained by the procedure {\sf MARK}. There exist a constant $\alpha\in(0,1)$ and $\beta>0$ such that
    \begin{equation*}
    \EnergyNorm{u- u_{\mesh_*}}^2 + \beta\, \eta^2_{\mesh_*}(u_{\mesh_*},\data)\le \alpha\left(\EnergyNorm{u - u_\mesh}^2 + \beta \,\eta^2_\mesh (u_\mesh, \data)\right),
    \end{equation*}
    where $\mesh_*$ is a refinement of the partition $\mesh$ obtained by applying the step {\sf REFINE}.
\end{thm}
\section{Numerical experiments}\label{sec:NumericalExperiments}
\begin{figure}[]
    \centering
    \subfloat[Maximum value of $\Lambda$ with respect to ${\texttt{Ndofs}}$ in the AVEM.\label{fig:Refinements:a}]{\includegraphics[width=.45\textwidth]{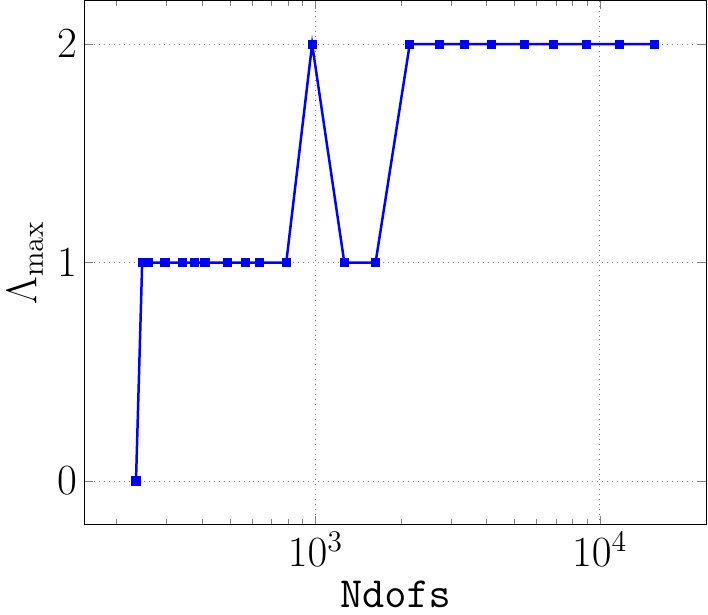}}\quad  
    \subfloat[The partition of the domain after 10 refinement iterations of the AVEM.\label{fig:Refinements:b}]{\includegraphics[width=.45\textwidth]{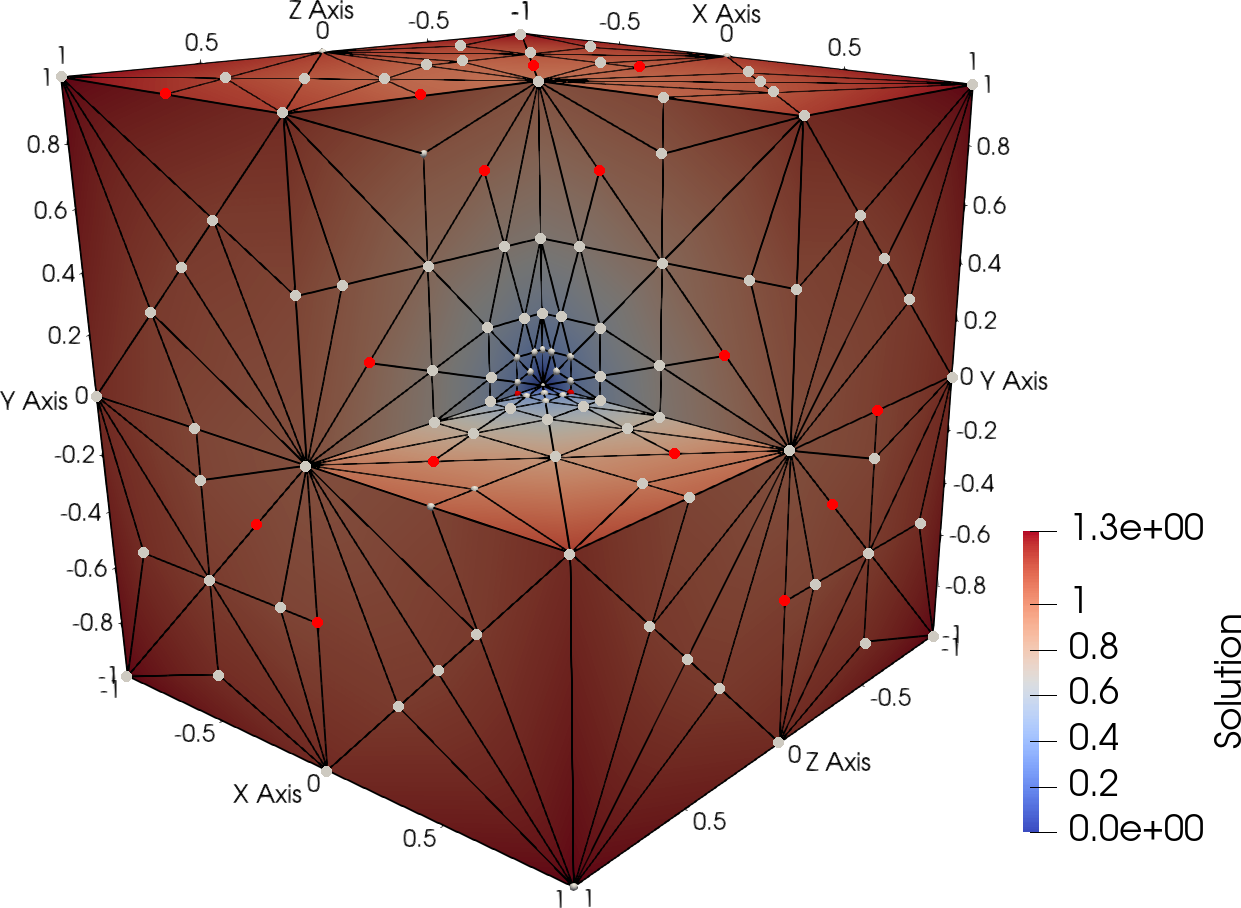}}\\
    \subfloat[Focus on the origin of the domain in the case of AFEM (1564 ${\tt Ndofs}$). All the nodes are proper nodes, and are depicted in white.\label{fig:Refinements:c}]{\includegraphics[width=.45\textwidth]{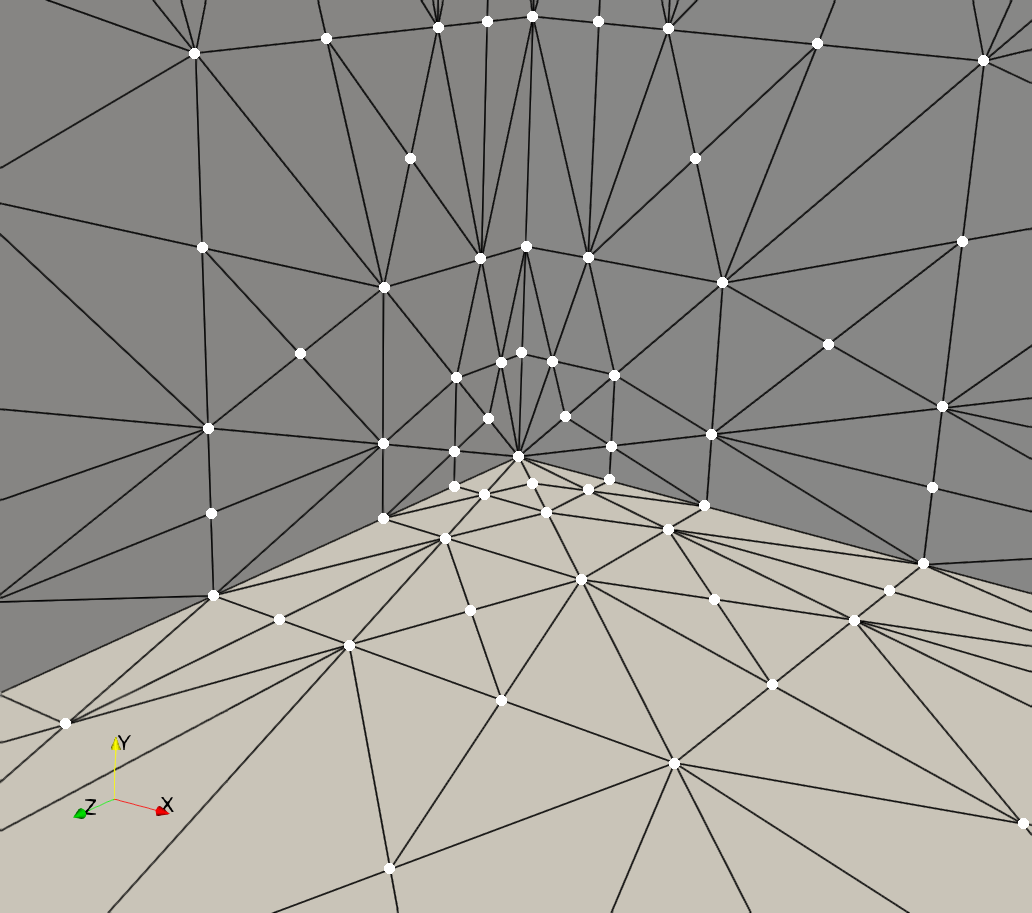}}\quad
    \subfloat[Focus on the origin of the domain in the case of AVEM (1628 ${\tt Ndofs}$). The proper nodes and the hanging nodes are depicted in white and red, respectively.\label{fig:Refinements:d}]{\includegraphics[width=.45\textwidth]{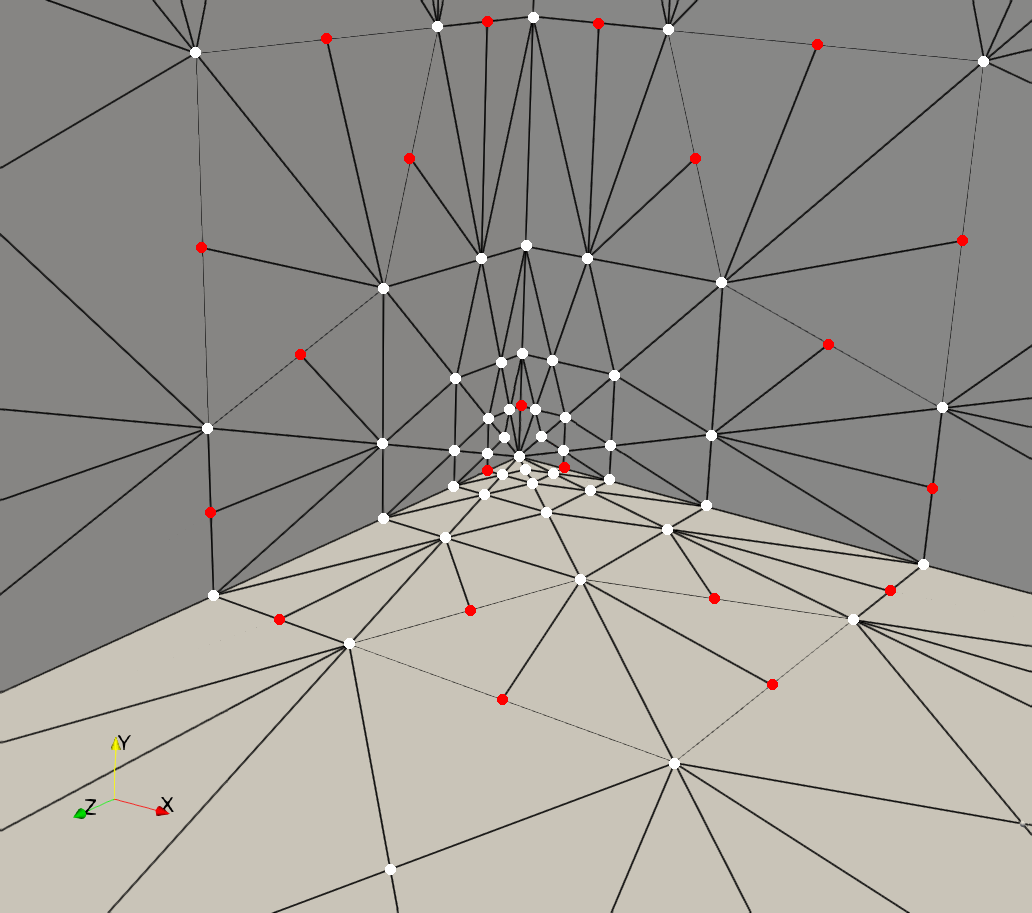}}
   
    \caption{The test solution and the refined meshes obtained.}
    \label{fig:Refinements}
\end{figure}
\begin{figure}[]
    \centering
    \subfloat[$H^1_{err}$  of { AFEM} and { AVEM} (dashed lines), the relative estimator $\eta_\mesh$ of { AFEM} and { AVEM} (solid lines) with respect to ${\texttt{Ndofs}}$.\label{fig:convergenceNdofs2:a}]{\includegraphics[width=.45\textwidth]{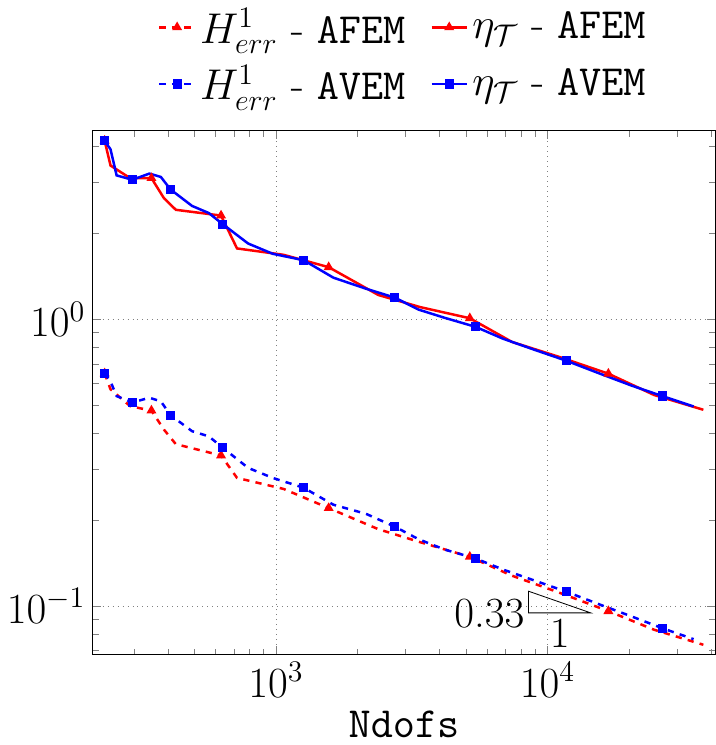}}\quad  
    \subfloat[The relative estimator $\eta_\mesh$ of { AVEM} (solid blue) and the stability term $S_\mesh$ (dashed green) with respect to ${\texttt{Ndofs}}$.\label{fig:convergenceNdofs2:b}]{\includegraphics[width=.45\textwidth]{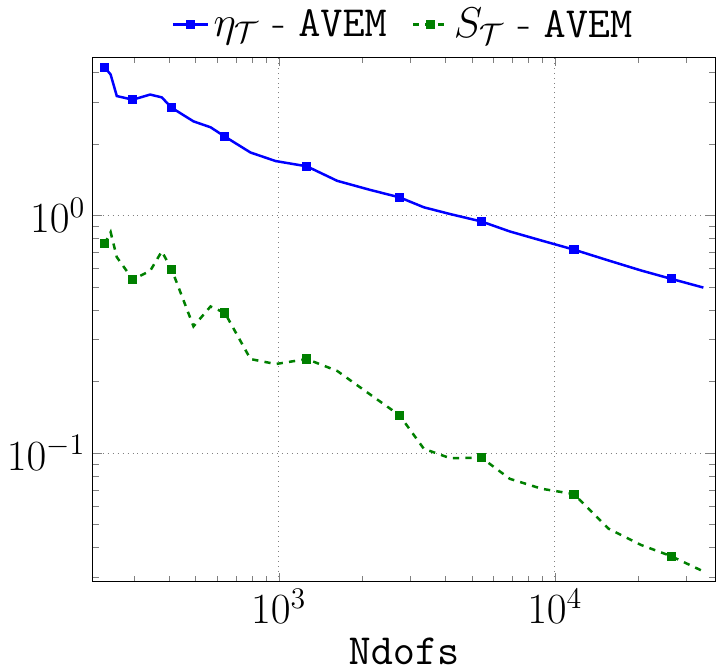}}\\
    \subfloat[${\texttt{Ndofs}}$ of { AFEM} and { AVEM} (dashed lines), number of mesh cells ($\# \mathcal{T}$) of { AFEM} and { AVEM} (solid lines) with respect to $H^1_{err}$.\label{fig:convergenceNdofs2:c}]{\includegraphics[width=.45\textwidth]{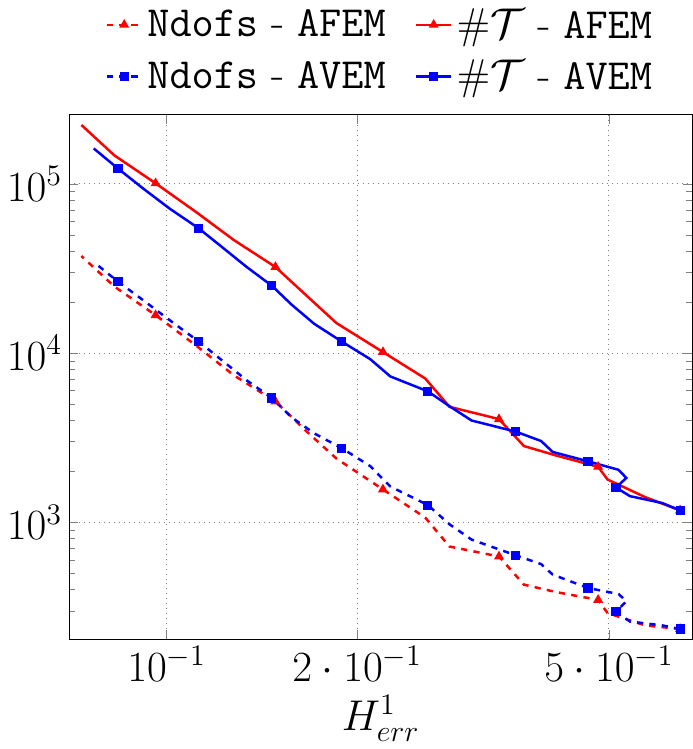}}\quad
    \subfloat[Number of marked cells ($\#\mathcal{T}_{mrk}$) in { AFEM} (dashed line) and number of refined cells ($\#\mathcal{T}_{ref}$) in { AFEM} and { AVEM} (solid lines) with respect to ${\texttt{Ndofs}}$.\label{fig:convergenceNdofs2:d}]{\includegraphics[width=.45\textwidth]{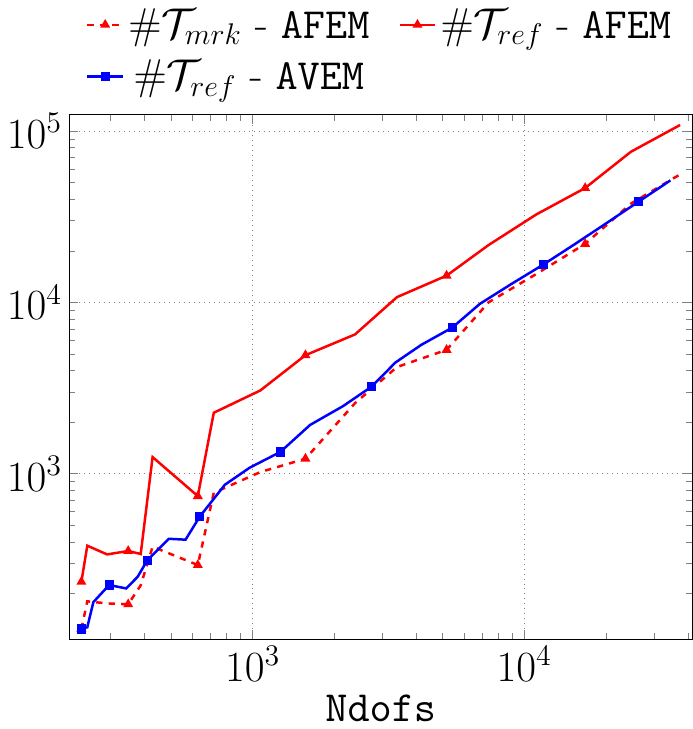}}
   
    \caption{Error analysis of the test.}
    \label{fig:convergenceNdofs2}
\end{figure}

We consider the Fichera test \cite{FicheraCornerDomain1, FicheraCornerDomain2} consisting in a three-dimensional domain $\Omega = (-1,1)^3\setminus (0,1)^3$, which is a natural extension of the bi-dimensional L-shape domain. We solve the Poisson problem \eqref{Initial_Problem} with $\K =  I$, where $I$ is the identity matrix, $c=1$, and Dirichlet boundary conditions such that the exact solution results
\begin{equation*}
    u_{ex}(\rho, \theta, \psi)= \rho^{\alpha}, 
\end{equation*}
with $\rho$, $\theta$ and $\psi$ the polar coordinates, and $\alpha\in(0,1)$. 

It is possible to take $u_{ex} \in W^k_p(\Omega)$, $\Omega\subset \mathbb{R}^d$, by exploiting the bound $k- \frac{d}{p} +\frac{d}{q}>0$ \cite[Theorem 6.29]{Bonito2024adaptive}, with $q=2$ and $p\ge 1$. In particular, there exists $p\in\left(\frac{d}{k+d/2},\frac{d}{k-\alpha}\right)$ such that the adaptive algorithm can achieve the optimal order of convergence, i.e. ${\texttt{Ndofs}}^{k/d}$, where $\texttt{Ndofs}$ indicates the number of degrees of freedom. In our case, we take $\alpha =1/2$, $d=3$, and $k=1$, and we get $p\in\left(\frac{6}{5},6\right)$. We recall that a VEM function is not known explicitly inside the elements $E$. Then, to compute the error, we project the solution onto $\mathbb{P}_1(E)$ via $\Pi^\nabla_E$, and we define the error as   
\begin{equation*}
     H^1_ {err}\coloneqq \frac{\left(\sum_{E\in \mesh}\normLE{\nabla(u_{\text{ex}}- \Pi^\nabla_E u_\mesh )}^2\right)^{1/2}}{\norm{\nabla u_{\text{ex}}}_{0, \Omega}}.
\end{equation*}
We set a maximum number of $\texttt{Ndofs}$ equal to $38000$ and the D\"{o}rfler parameter $\theta$ for the {\sf MARK} step equal to $0.5$. We carry on a study on the value $\Lambda$. In particular, we analyze the two extreme cases in which $\Lambda=0$ and $\Lambda=10$. We recall from Remark~\ref{Remark:LambdaMax} that setting $\Lambda= 0$ we force our mesh not to admit hanging nodes, and the Adaptive Virtual Element Method algorithm converges to the classical Adaptive Finite Element Method built on tetrahedral elements. On the other hand, the value $\Lambda =10$ does not represent a limit in the number of hanging nodes in the mesh for this test. Indeed, it can be seen from the simulation that the maximum global index of the partition floats between the values of one and two, see Figure~\ref{fig:Refinements:a}.  
We label AFEM and AVEM the cases $\Lambda=0$ and $\Lambda=10$, respectively.

The solution is shown in Figure~\ref{fig:Refinements:b}, with a focus on the resulting refined meshes in the origin of the domain reported in Figures~\ref{fig:Refinements:c}-\ref{fig:Refinements:d}. 
As shown in Figure \ref{fig:Refinements:d} the number of hanging nodes grows in the proximity of the origin, where the solution is less regular. This is the region where the AFEM propagates more the refinements (see Figure \ref{fig:Refinements:c}) and the number of cells of the AFEM is greater with respect to the AVEM.

Figure~\ref{fig:convergenceNdofs2:a} reports the convergence plots of  $H^1_{err}$ and $\eta_\mesh$ with the optimal decay rate for both the AFEM and AVEM.  
Figure~\ref{fig:convergenceNdofs2:b} represents a numerical confirmation of the bound found in Proposition \ref{Bound}, in which the stability term $S_\mesh(u_\mesh,u_\mesh)$ results bounded by $\eta_\mesh(u_\mesh,\data)$.
The quantities $S_\mesh$ and $\eta_\mesh$ represented in all the plots are scaled with $\norm{\nabla u_{\text{ex}}}_{0, \Omega}$.

We observe that the AVEM generates a partition with the 27,33\% less of three-dimensional cells $\#\mesh$ (161937 vs 222866) with respect to AFEM, reaching the same $ H^1_{err}$ order of magnitude, see Figure \ref{fig:convergenceNdofs2:c}. 
Finally, Figure~\ref{fig:convergenceNdofs2:d} shows that, despite the  AFEM and the  AVEM have approximately the same number of marked three-dimensional cells, having fixed the ${\tt Ndofs}$, the number of refined elements in the AVEM is around 30\% less than the AFEM. We highlight that this also implies a proportional decrease in the computational cost related to the refinement procedure.

\section{Conclusions} 
In this paper, we investigate the three-dimensional Adaptive Virtual Element Method (AVEM) algorithm, supported by computational experiments. 
We focus on tetrahedra with hanging nodes, i.e., aligned faces and edges, extending the results obtained for the two-dimensional case \cite{AVEMstabfree}.

By introducing a new proof of the scaled Poincaré inequality, we provide an extension to the three-dimensional case of the stabilization-free a posteriori bound for the energy error and compare the proposed algorithm with the Adaptive Finite Element Method (AFEM).
The AVEM shows a reduced number of mesh elements in the refinement procedure, for a fixed number of degrees of freedom and error, leading to a proportional decrease in the computational cost associated with the refinement process.

\section*{Acknowledgments} The authors wish to thank Professor Claudio Canuto for the enlightening discussions on adaptivity. 

SB acknowledges the funding by the European Union through project Next Generation EU, M4C2, PRIN 2022 PNRR project P2022BH5CB\_001 PG-stab ``Polyhedral Galerkin methods for engineering applications to improve disaster risk forecast and management: stabilization-free operator-preserving methods and optimal stabilization methods'' and the
MUR PRIN project 20204LN5N5\_003 AdPolyMP - ``Advanced polyhedral discretisations of heterogeneous PDEs for multiphysics problems''.

DF performed this research in the framework of the Italian MIUR Award ``Dipartimenti di Eccellenza 2018-2022'' granted to the Department of Mathematical Sciences, Politecnico di Torino CUP\_E11G18000350001.

FV has been partially funded by the PRIN project ``FREYA - Fault REactivation: a hYbrid numerical Approach'', Project ID 2022MBY5JM CUP D53D23005890006.

DF and FV thank the INdAM - GNCS Project CUP\_E53C23001670001.

All the authors are members of the Italian INdAM-GNCS research group.
\bibliographystyle{acm}
\bibliography{references}

\end{document}